\documentclass[reqno]{amsart}
\usepackage{amsmath, amssymb}
\usepackage{url}
\usepackage{hyperref}
\usepackage{color}

\begin{document}

\title[Autonomous dynamics]
{On Partially Hyperbolic diffeomorphisms in dimension three\\ via a notion of autonomous dynamics
}

\author[]{Souheib Allout and Kambiz Moghaddamfar
}
      
\address{Souheib Allout \newline
Ruhr-Universitaet Bochum, Germany, and \newline Faculty of mathematics, USTHB, Algiers, Algeria}
\email{souheib.allout@rub.de}

\address{Kambiz Moghaddamfar \newline
UMPA, ENS de Lyon, France, and \newline Sharif university of technology, Tehran, Iran
 }
\email{kambiz.moghaddamfar@ens-lyon.fr}

\dedicatory{}

%\thanks{Submitted April 1, 2020. Published January 7, 2021.}
%\subjclass[2010]{35A15, 35B09, 35B50, 35D30}
\date{\today}

\begin{abstract}
We introduce a notion of autonomous dynamical systems and apply it to prove rigidity of partially hyperbolic diffeomorphisms on closed compact three-manifolds under some smoothness hypothesis of their associated framing. 
\end{abstract}

\maketitle
\numberwithin{equation}{section}
\newtheorem{theorem}{Theorem}[section]
\newtheorem{proposition}{Proposition}[section]
\newtheorem{lemma}[theorem]{Lemma}
\newtheorem{corollary}[theorem]{Corollary}
\newtheorem{remark}[theorem]{Remark}
\newtheorem{fact}[theorem]{Fact}
\allowdisplaybreaks

\def\C{\mathbb{C}}
\def\R{\mathbb{R}}
\def\Q{\mathbb{Q}}

\def\H{\mathbb{H}}
\def\Z{\mathbb{Z}}

\def\N{\mathbb{N}}

\def\F{\mathcal {F}}

\def\U{\sf{U}}
\def\SU{\sf{SU}}

%%%%%%%%%%%%%%%%%%%%%%%%%%%%%%%%%%%%%%%%%%

\def\P{\mathbb{P}}
\def\S{\mathbb{S}}
\def\D{\mathbb{D}}

\def\T{\mathbb{T}}
\def\M{\mathcal{M}}

%%%%%%%%%%%%%%%%%%%%%%%%%%%%%%%%%%%%%%%%%%

\def\CC{\mathcal{C}}

\def\la{\lambda}

\def\Iso{\sf{Iso}}
\def\Pro{\sf{Proj}}
\def\Aff{\sf{Aff}}
\def\Sim{\sf{Sim}}

\def\Euc{\mathsf{Euc}}

\def\Spa{\sf{Span}}

\def\Jac{\sf{Jac}\;}

\def\Met{\mathcal Met}

\def\mul{\sf{mul}}
\def\B{\mathcal{B}}

\def\I{\sf{\bf I}}

\def\SS{\mathcal S}

\def\E{\mathcal E}

\def\EE{\varepsilon}

%%%%%%%%%%%%%%%%%%%%%%%%%%%%%%%%%%%%%%%%%%%%%%%%%%%%%%%%%%%%%%%%%
%   Abbrev.
%%%%%%%%%%%%%%%%%%%%%%%%%%%%%%%%%%%%%%%%%%%%%%%%%%%%%%%%%%%%%%%%%

\def\Aut{{\sf{Aut}}}
\def\DF{{\sf{DF}}}
\def\Int{{\sf{Int}}}

\def\MCG{{\sf{MCG}}}
\def\PGL{{\sf{PGL}}}
\def\PSL{{\sf{PSL}}}
\def\GL{{\sf{GL}}}
\def\O{{\sf{O}}}
\def\SO{{\sf{SO}}}
\def\SU{{\sf{SU}}} 
\def\SL{{\sf{SL}}}
\def\Sp{{\sf{Sp}}}
\def\End{{\sf{End}}}
\def\Mat{{\sf{Mat}} }
\def\Iso{{\sf{Iso}} }
\def\Stab{{\sf{Stab_1}} }
\def\Exo{{\sf{Exercise}} }
\def\Lie{{\sf{Lie}} }
\def\Or{{\sf{Or}}}
\def\Span{{\sf{Span}} }
\def\Ad{{\sf{Ad}} }

\def\Mink{{\sf{Mink}} }
\def\dS{{\sf{dS}} }
\def\AdS{{\sf{AdS}} }
 
\def\tr{{\sf{tr}}}
\def\Tr{{\sf{Tr}}}
\def\rk{{\sf{rk}}}
\def\det{{\sf{det}}}

\def\p{{\mathfrak{p}}}
\def\k{{\mathfrak{k}}}
\def\t{{\mathfrak{t}}}
\def\a{{\mathfrak{a}}}
\def\g{{\mathfrak{g}}}
\def\sll{{\mathfrak{sl}}}
\def\e{{\mathfrak{e}}}
\def\so{{\mathfrak{so}}}
\def\sp{{\mathfrak{sp}}}
\def\n{{\mathfrak{n}}}
\def\s{{\mathfrak{s}}}
\def\h{{\mathfrak{h}}}
\def\i{{\mathfrak{i}}}
\def\z{{\mathfrak{z}}}
\def\z{{\mathfrak{z}}}
 \def\q{{\mathfrak{q}}}
\def\I{{\mathfrak{I}}}

 \newcommand{\ghani}[1]{\textcolor{red}{\sf ghani: #1}}
\newcommand{\Souheib}[1]{\textcolor{blue}{\sf  X: #1}}

\tableofcontents
\section{\textbf{Introduction and statement of main results}}
The study of partially hyperbolic diffeomorphisms in dimension three, both in the flexible and rigid settings, has been widely explored and it is still an active area of research. Rigidity of these diffeomorphisms, satisfying some additional requirements, has been investigated in both smooth and topological frameworks. In this paper, which is motivated by \cite{CPH20}  (and also \cite{MM}), we provide smooth "rigid" classifications of partially hyperbolic diffeomorphisms satisfying some additional hypothesis. We show then how to extend this to the $C^1-$regularity case. 
\subsection{An Introductory example}\label{1.1} (see more details in Section \ref{C0}) The better is perhaps to start with the following construction illustrating the content of the present paper.
Let $A:\T^2\to \T^2$ be a hyperbolic toral linear automorphism. Let $c: \T^2 \to \S^1$ be a smooth map and consider the \textit{skew product} $\varphi$ on $M=\T^2\times \S^1$ defined by: $$\varphi(x, \theta)=(A(x), \theta + c(x))$$
One can construct continuous vector fields $X_s, X_u,$ and $X_c$ on $M$ satisfying $\varphi_*(X_s)=\lambda_s X_s, \varphi_*(X_u)=\lambda_u X_u$ and $\varphi_*(X_c)=\lambda_c X_c$. The vector field $X_c$ corresponds to the $\S^1$ factor, while $X_s$ and $X_u$ correspond to linear vector fields on $\T^2$ and $0<\lambda_s<1<\lambda_u$ are the eigenvalues of $A$. The existence of $X_s$ and $X_u$ follows from basics of partially hyperbolic dynamics theory which will be the central subject in this paper.\\
Our map $\varphi$ is a special partially hyperbolic diffeomorphism in the sense that it has continuous, even constant, Lyapunov exponents (and splittings)! This situation corresponds to our second topic here: autonomous dynamical systems, to mean there exists a framing where the derivative cocycle map is constant. We prove a differentiable rigidity of partially hyperbolic autonomous systems in dimension three as well as autonomous systems in dimension two.\\
To get a flavour of the results and methods in this paper, let us assume that $X_s$ and $X_u$ are $C^1$. Their bracket $Z=[X_s,X_u]$ is a well defined $C^0$ vector field. One observes that $Z= a X_c$ for a $\varphi-$invariant function $a$. But, $a$ being $\S^1-$invariant implies that it is in fact constant. We also have $[X_s, X_c] = [X_u, X_c] = 0$, and 
therefore, $\{X_s, X_c, X_u\}$ generates a three dimensional Lie algebra isomorphic to the Heisenberg algebra if $a\neq 0$ and to $\R^3$ otherwise. Thus, $M$ is a quotient of the simply connected Lie group, associated to this algebra, by a lattice. By algebraic topological arguments, one excludes the Heisenberg case and, thus, $\R^3$ acts transitively and locally freely on $M$. One, then, shows that $\varphi$ is $C^1-$conjugate to an affine partially hyperbolic automorphism on $\T^3$. Its linear part (as an element of $\GL(3,\Z)$) is, up to a finite cover, $C^1-$conjugate to a direct product: $$(x,\theta)\in \T^2\times \S^1\to (A(x), \theta)\in \T^2\times\S^1$$
Obviously, this is a special situation that happens for rare maps $c:\T^2\to \S^1$; even if the $C^0-$vector fields ($X_s, X_u$ and $X_c$) exist for all.\\\\
Let us emphasize that our regularity assumption (for vector fields) is merely $C^1$, and not $C^2$ as it is usually the case for structures involving a kind of curvature...! For this it is worthwhile to recall that a $C^1$ vector field generates a $C^1$ flow (not just a $C^0$ one). Also, applying any definition of brackets (for instance as commutators of operators on $C^{\infty}$ functions), the bracket of two $C^1$ vector fields is a well defined $C^0$ vector field.\\ 
Actually we think that the previous outlined rigidity extends to the case where the vector fields are Lipschitz?! Indeed, Frobenius' theorem is valid with this regularity (\cite{sim}), and therefore also is the theorem of Palais on integration of Lie algebra actions since its proof is based on Frobenius (see the appendix \ref{appendix} for more details).
\subsection{} After this warming up example,
let us give more precise definitions.
Let $M$ be a compact smooth manifold of dimension $n$ endowed with a parallelization (framing) $\mathcal{F}$ of
its tangent bundle, $i.e.$ a system of vector fields defining a basis of each tangent space and let $\varphi$ be a $C^1$ diffeomorphism of $M$. The derivative cocycle is the map: $$x\in M\to C_{\varphi}(x)\in \GL(n,\R)$$ where $C_{\varphi}(x)$ is the matrix of the derivative $D_x\varphi: T_xM\to T_{\varphi(x)}M$, when these linear spaces are endowed with the bases $\mathcal{F}(x)$ and $\mathcal{F}(\varphi(x))$ respectively.\\\\
We say that $\varphi$ is \textit{autonomous} with respect to $\mathcal{F}$ if the cocycle map $C_{\varphi}$ is constant. The same definition applies to a $C^1$ action of a group $G$ to mean that any element of it is autonomous with respect to $\mathcal{F}$.\\\\
Now, we say that a $G-$action on $M$ is autonomous if it is autonomous with respect to some parallelization $\mathcal{F}$ on $M$. In this case we have a representation $G\to \sf{GL}(n,\R)$ associating to each $g\in G$ the cocycle matrix of $g$ (acting on $M$). This is essentially inspired from \cite{Z1} and \cite{Z2} where the terminology \textit{autonomous} was introduced. Our contribution here is to prove results towards the classification of such autonomous systems in some cases. As one may expect, the regularity of the framing is relevant and it is natural to ask when is possible to have constant matrix? $i.e$ when is it possible for a diffeomorphism to be autonomous?
\subsubsection{Regularity of the framing} We will always assume the manifold $M$ and the diffeomorphism $\varphi$ to be smooth. Now, our rigidity results require the framing to be at least of class $C^1$. Continuous framings are also interesting, but we will not deal with them in the present article. The introductory example in Section \ref{1.1} shows how abundant they are. Less than $C^0$, say measurable framings, seem to be not restrictive.
\subsection{Examples} Here we give some examples of autonomous actions that will be explained in more details in the next section:
\subsubsection{\textbf{Toral affine automorphisms:}}
Consider the space $\T^n=\R^n/\Z^n$. The group $\GL(n,\Z) \ltimes \T^n$ acts naturally on $\T^n$ by affine automorphisms (where $\GL(n,\Z)$ is the group of integer matrices with determinant $\pm 1$). This action is autonomous with respect to the canonical framing on $\T^n$. The cocycle matrix of each element $A\in \GL(n,\Z) \ltimes \T^n$ coincides with matrix of the linear part of $A$.
\subsubsection{\textbf{The flat case:}} Let $\mathcal{F}$ be a $C^{k}, \ k\geq 1$ framing given by linearly independent vector fields $X_1,X_2, \ldots, X_n$ of a compact smooth manifold of dimension $n$. We say that $\mathcal{F}$ is flat if $[X_i,X_j]=0$ for every $i,j\in \{1,2,\ldots,n\}$. In this case the flows $\phi_1^t, \phi_2^t,\ldots,\phi_n^t$ (of the vector fields $X_1,X_2, \ldots, X_n$ respectively) commute and generate a transitive locally free $C^{k}$ action of $\R^n$ on $M$. Hence $M\approx \R^n/\Gamma$ where $\Gamma$ is a discrete subgroup with compact quotient. This implies that $M\approx\T^n$ and we have:
\begin{proposition} \label{1.1} Let $M$ be a compact manifold of dimension $n$ with a flat $C^{k}, \ k\geq 1$ framing $\mathcal{F}$. Assume that $\varphi$ is a smooth diffeomorphism of $M$ which is autonomous with respect to $\mathcal{F}$. Then $\varphi$ is $C^k$ conjugate to an affine automorphism of the $n-$dimensional torus $\T^n$.
\end{proposition}
Proof: We have that $[X_i,X_j]=0$ for every $i,j\in \{1,2,\ldots,n\}$. This implies that the flows $\phi_1^t, \phi_2^t,\ldots,\phi_n^t$ generated by these vector fields define all together a $C^k$ action of $\R^n$ on $M$ which is locally free because the vector fields define a basis at each point and it is transitive because the orbits in this case are open and $M$ is connected. Hence $M=\R^n/\Gamma$ where $\Gamma$ is a discrete subgroup and the framing $\mathcal{F}$ is induced by projecting a parallel framing $\mathcal{P}$ of $\R^n$. The diffeomorphism $\varphi$ has a lifting $\widetilde{\varphi}:\R^n\to \R^n$ satisfying $\varphi\circ \pi=\pi \circ \widetilde{\varphi}$ where $\pi:\R^n\to\R^n/\Gamma$ is the natural projection. Hence $\widetilde{\varphi}$ is autonomous with respect to $\mathcal{P}$ with cocycle matrix $M_{\varphi}=M_{\widetilde{\varphi}}$. This implies that $\widetilde{\varphi}$ is affine. \begin{flushright}
				$\square$
			\end{flushright}
\subsubsection{\textbf{The algebraic case:}}
Let $G$ be a Lie group and consider $M = G/\Gamma$ where $\Gamma$ is a discrete subgroup acting by right translations on $G$. Take a basis $\mathcal{B}$ of the Lie algebra $\sf{Lie}(G)$ and right translate it to get a framing of $G$. The projection of this framing defines a framing $\mathcal{F}$ of $G/\Gamma$. Then the left $G-$action on $G/\Gamma$ is autonomous. In fact, for $g\in G$, its derivative cocycle when acting on $G/\Gamma$ is just the matrix of $\sf{Ad}(g)$ acting on $\sf{Lie}(G)$ with respect to the basis $\mathcal{B}$.\\\\
Let now $A:G\to G$ be an automorphism of $G$ such that $A(\Gamma)=\Gamma$ and $g\in G$. Then $L_g\circ A$ induces a diffeomorphism on $G/\Gamma$ which is autonomous with respect to $\mathcal{F}$ and having cocycle matrix $M_{\sf{Ad(g)}}M_A$ where $M_{\sf{Ad}(g)}$ and $M_A$ are the cocycle matrices of $\sf{Ad}(g)$ and $A$ respectively. We call a composition like $L_g\circ A$ an \textit{affine automorphism}. We will give more details about all this.
\subsection{Partially hyperbolic diffeomorphisms}Let $M$ be a compact smooth manifold of dimension $n$ and let $f:M\to M$ be a diffeomorphism of $M$. Then $f$ is said to partially hyperbolic if the tangent bundle $TM$ splits continuously into the direct sum into three $f-$invariant sub-bundles $E_s$, $E_c$, and $E_u$ such that $TM=E_s\oplus E_c\oplus E_u$ and there exists a Riemannian metric on $M$ such that for all unitary vectors $v_s, v_c, v_u$ in $E_s, E_c$, and $E_u$ respectively, we have:
\begin{itemize}
 \item $\left \| Df(v_s) \right \|<1, \ \text{and} \ \left \| Df(v_u) \right \|>1$
\item $\left \| Df(v_s) \right \|<\left \| Df(v_c) \right \|<\left \| Df(v_u) \right \|$ 

\end{itemize}

\subsubsection{\textbf{A rigidity hypothesis:}} In the paper \cite{CPH20}, which is recently presented to classify partially hyperbolic diffeomorphisms in dimension three, the authors considered partially hyperbolic diffeomorphisms $\varphi:M\to M$ satisfying a rigidity hypothesis, $i.e$, having constant Lyapunov exponents. It turns out that this is equivalent to the fact that the diffeomorphism $\varphi$ is autonomous with respect to a framing $\mathcal{F}$ given by three vector fields $X_s$, $X_c$, and $X_u$ such that $\varphi_*(X_s)=\lambda_s X_s$, $\varphi_*(X_c)=\lambda_c X_c$, and $\varphi_*(X_u)=\lambda_u X_u$ where $\left|\lambda_s \right|<1$, \ $\left|\lambda_u \right|>1$, and $\left|\lambda_s \right|<\left|\lambda_c \right|<\left|\lambda_u \right|$.\\\\
In \cite{CPH20} they provide a classification with the assumption that the vector fields $X_s, X_c,$ and $X_u$ are smooth, and that the diffeomorphism $\varphi$ is either transitive or real analytic. We give here a more explicit classification assuming only $C^1$ regularity without the transitivity (or analyticity) assumption.\\\\
Let $M$ be a compact connected three manifold and $\varphi:M\to M$ be a partially hyperbolic diffeomorphism satisfying the rigidity hypothesis $i.e$ $\varphi$ is autonomous with respect to a framing $\mathcal{F}$ given by three vector fields $X_s, X_c$, and $X_u$ such that the cocycle matrix of $\varphi$ equals $\begin{pmatrix}
\lambda_s &0 &0 \\ 
0 &\lambda_c &0 \\
0 &0 &\lambda_u
\end{pmatrix}$ where $\left|\lambda_s \right|<1$, \ $\left|\lambda_u \right|>1$, and $\left|\lambda_s \right|\leq\left|\lambda_c \right|\leq\left|\lambda_u \right|$. In this case we have:
\begin{theorem}
\label{3D} Assume that the vector fields $X_s, X_c$, and $X_u$ are $C^1$. Then, up to finite power and cover, $\varphi$ is $C^1-$conjugate to an affine automorphism $L_g\circ A:G/\Gamma\to G/\Gamma$ where $A\in \sf{Aut}(G)$ satisfying $A(\Gamma)=\Gamma$, $L_g$ is a left translation, and $\Gamma\subset G$ is a cocompact lattice of the simply connected three dimensional Lie group $G$. In addition, $G$ is one of the four Lie groups: \ $\R^3$, the Heisenberg group (denoted $\sf{Heis}_3$), $\sf{Lor}(1,1)=\SO(1,1)\ltimes \R^2$, or \ $\widetilde{\sf{SL}}(2,\R)$.

\end{theorem}
See Sections  \ref{5.6} and \ref{5.2} for more details. We will present first an elementary proof (relying somehow on Sard's theorem) assuming $C^{\infty}$ regularity of the framing, showing that the conjugacy is smooth. The point is to emphasize that the $C^1$ regularity is by no means trivial which requires a different proof. After that, in Section \ref{C0},   we show how to construct counterexamples to Theorem \ref{3D} in the $C^0-$regularity case as quickly mentioned  Section \ref{1.1}   (see Corollary \ref{6.2}).

\subsubsection{\textbf{Smooth vs topological classification, related results:}}

${}$

$\bullet$ One observation is that, in our case where $\lambda_c=1$ and $\lambda_s$, $\lambda_u$ have the same sign, the diffeomorphism $\varphi:M\to M$ preserves a Lorentzian metric on $M$ given by the quadratic form $g=dx_s dx_u + dx_c^2$ with respect the framing generated by the vector fields $X_s, X_c, X_u$. This implies, in particular, that the isometry group of $(M,g)$ is non-compact. In the paper \cite{CF},  C. Frances gives a topological (and geometrical) classification of closed Lorentzian three-manifolds with non-compact isometry groups. He showed that such a manifold must be homeomorphic to one of the following, assuming at least $C^9$ regularity:
\begin{itemize}
\item A quotient $\widetilde{\SL}(2,\R)/\Gamma$ where $\Gamma\subset \widetilde{\SL}(2,\R)$ is a cocompact lattice.
\item The torus $\T^3$, or, a suspension of $\T^2$ by a parabolic or a hyperbolic linear automorphism.
\end{itemize}
Theorem \ref{3D} classifies (assuming only $C^1$ regularity) not only the manifold $M$, but even the dynamics for being algebraic. In particular, comparing to C. Frances' results, the Lorentzian metric $g$ on $M$ (constructed above) is in fact locally homogeneous.\\
$\bullet$ Another remark is that in our situation the case $\sf{Heis}_3/\Gamma$ corresponds to a suspension of $\T^2$ by a parabolic element, and the case $\sf{Lor}(1,1)/\Gamma$ corresponds to a suspension of $\T^2$ by a hyperbolic element, in the results of C. Frances.\\\\
$\bullet$ In \cite{MM},  Mion-Mouton proves a rigidity of partially hyperbolic diffeomorphisms in another direction. He classifies the three-dimensional partially hyperbolic diffeomorphisms whose stable, unstable, and central distributions $E_s, E_u,$ and $E_c$
are smooth, such that $E_s\oplus E_u$ is a contact distribution, and whose non-wandering set equals the whole manifold. He shows that up to a finite quotient or a finite power,
they are smoothly conjugate either to the time-one map of an algebraic contact Anosov flow, or to an affine partially hyperbolic automorphism of a nil-manifold.\\\\
$\bullet$  In \cite{BGP} and \cite{BKP},  the authors construct new families of partially hyperbolic diffeomorphisms in dimension three using topological (surgery) methods. This shows the flexible side of this theory. Whereas, for example in \cite{BW05} and \cite{BZ},  the authors provide some rigid topological classifications, for instance in \cite{BZ} they classify  partially hyperbolic diffeomorphisms assuming transitivity and the existence of $1-$dimensional neutral central direction. In this paper, we give, particularly, \textit{smooth} classifications. \\\\
$\bullet$  Lastly, let us say  that the rigidity hypothesis (constancy of Lyapunov exponents), which is also assumed in \cite{CPH20},   seems at first to be restrictive,   but it is natural  enough  because  this hypothesis is,  in particular, 
true for algebraic systems.  In fact, our framework can also apply to handle higher dimension situations. In opposite,  as shown in \cite{MM}, a merely smoothness hypothesis of distributions, is far from leading to prove rigidity.   

%In fact, pour framework can apply to also handle higher dimensions situations. In opposite,  as shown in \cite{MM}, a merely smoothness hypothesis of %distributions, is far from leading to prove rigidity.   

%the case of algebraic to consider this case (especially when one wants to provide an algebraic classification), and yields a %nice highly non-trivial classification as shown in \cite{CPH20}.
%$\bullet$  Lastly, let us mention that the rigidity hypothesis (constancy of Lyapunov exponents), which is also assumed in %\cite{CPH20}, seems at first to be restrictive, but it is natural to consider this case (especially when one wants to provide an %algebraic classification), and yields a nice highly non-trivial classification as shown in \cite{CPH20}.

\subsection{Autonomous diffeomorphisms in dimension two}
Let $S$ be a compact connected smooth surface and $\mathcal{F}$ a $C^{\infty}$ framing of $S$. Then $S$ is necessarily diffeomorphic to the torus $\T^2$. Let $\varphi:S\to S$ be an orientation preserving diffeomorphism which is autonomous with respect to $\mathcal{F}$ and having cocycle matrix $M_{\varphi}$. Then:
\begin{theorem}
\label{2D} We have the following classification:
\begin{itemize}
\item If the cocycle matrix $M_{\varphi}$ is hyperbolic, then $\varphi$ is smoothly conjugate to a linear Anosov automorphism of $\T^2$.
\item If it is elliptic, then $\varphi$ is smoothly conjugate to an isometry of a flat torus.
\item If it is parabolic (with eigenvalues equal to $1$), then either $\varphi$ is smoothly conjugate to an \textbf{affine} parabolic automorphism or, a finite power of $\varphi$ is smoothly conjugate to a \textbf{linear} parabolic automorphism of $\T^2$.
\end{itemize} 

\end{theorem}

The hyperbolic and the elliptic cases are rather straightforward. The parabolic case is the challenging one. More facts and details are in the proof section, in particular, we explain the two possibilities (linear, affine) in the parabolic case.
\begin{remark}\label{regularity}(About regularity) In this two dimensional case we assumed the framing to be smooth. Observe, however  that exactly as in the three-dimensional partially hyperbolic case, $C^1$ regularity is enough in the two-dimensional hyperbolic case. Actually, it seems to us that the same rigidity for Anosov autonomous systems holds  with just a $C^0$ assumption of the framing (and in fact for any dimension). Also the classification in the elliptic case (dimension two) seems to be still true in the $C^0$ case (the idea would be to apply the $C^0$ uniformization theorem). Lastly, the parabolic case requires $C^3$ regularity but we don't  know if it is optimal.
    
\end{remark}

\subsection{Questions} One may derive some questions:

-  In Theorem \ref{3D}, what happens if the vector 
fields $X_s$, $X_c$, and $X_u$ are only Lipschitz?

- In Theorem \ref{2D}, we can assume that the surface $S$ is not closed. In this case, if $S$ has finite volume (given by the framing $\mathcal{F}$) and the diffeomorphism $\varphi$ has hyperbolic cocycle matrix, then the framing $\mathcal{F}$ is flat defining a translation structure on $S$. Suppose that $S$ has finite topology (of finite type) then, is it isomorphic, as a surface with a translation structure, to a translation surface with a finite number of points removed (including the singularities)?  

- As suggested by the referee, it is natural and interesting to  also study autonomous diffeomorphisms
on $3-$manifolds which  are not partially hyperbolic. We hope to come back to this topic in a forthcoming work, but 
we can already observe that 
 the unipotent case in dimension $2$ was  quite  laborious.

%We would like to express our sincere gratitude to Professor Ghani Zeghib for proposing this project to us, his guidance, and his inspirational ideas.

\section{\textbf{General constructions}}
In this section we will give some basic constructions and definitions, mainly about autonomous systems.

\subsection{Some properties}\label{properties} Let $M$ be a compact manifold of dimension $n$ with framing $\mathcal{F}$ and $\varphi$ be an autonomous diffeomorphism with respect to $\mathcal{F}$ having cocycle matrix $M_{\varphi}$. Let $T$ be a tensor field on $M$ of a given type (for example a Riemannian metric, a differential form, a vector field, \ldots) and consider the map $\widetilde{T}:M\to \mathcal{T}$ where $$\mathcal{T}=\{\text{space of tensors on} \ \R^n \ \text{having the same type as} \ T\}$$ such that $\widetilde{T}$ associates to $z\in M$ the components of $T(z)$ with respect to the framing $\mathcal{F}(z)$. We know that $\varphi$ acts naturally on $T$ and $M_{\varphi}$ acts on $\mathcal{T}$, say via a representation $\rho_{\mathcal{T}}:\GL(n,\R)\to \GL(\mathcal{T})$. The map $\widetilde{T}$ is equivariant, that is: $$\widetilde{T}\circ \varphi=\rho_{\mathcal{T}}(M_{\varphi})\widetilde{T}$$
Let $g_0$ be the usual inner product in $\R^n$ (with respect to the canonical basis) and consider the tensor $T_0$ on $M$ such that the map $\widetilde{T}_0:M\to\sf{Sym}^{+}$ is constant and equals the identity matrix (the matrix of $g_0$ in the canonical basis). This defines a Riemannian metric on $M$ for which the basis $\mathcal{F}(z)$ is orthonormal for every $z\in M$. Consider now the associated Riemannian volume form $\omega_0$ on $M$. We have that $\varphi_*\omega_0=\det(M_{\varphi})\omega_0$. On the other side we have: $$\int_M\varphi_*\omega_0=\int_M \det(M_{\varphi})\omega_0=\pm\int_M \omega_0$$
This implies that $\det(M_{\varphi})=\pm 1 $. In particular, if $\varphi$ preserves the orientation (given by $\mathcal{F}$) then $M_{\varphi}\in \SL(n,\R)$.\\\\
Consider a basis $\mathcal{B}$ of $\R^n$ and let $A_{\mathcal{B}}$ be the matrix that sends the canonical basis to the basis $\mathcal{B}$. Define the framing $\mathcal{F}_{\mathcal{B}}$ such that $A_{\mathcal{B}}$ is the matrix that sends the basis $\mathcal{F}(x)$ of $T_xM$ to the basis $\mathcal{F}_{\mathcal{B}}(x)$ for every $x\in M$. Then $\varphi$ is autonomous with respect to $\mathcal{F}_{\mathcal{B}}$ with cocycle matrix $A_{\mathcal{B}}\circ M_{\varphi}\circ A_{\mathcal{B}}^{-1}$.
\subsection{Lie groups} Let $G$ be a simply connected Lie group of dimension $n$ with Lie algebra denoted by $\sf{Lie}(G)$, and let $\Gamma\subset G$ be a discrete subgroup. The quotient space $G/\Gamma$ is a homogeneous manifold of dimension $n$ whose points are the classes $x\Gamma$, for $x\in G$.
The group $G$ acts naturally on the left on the quotient $G/\Gamma$ as follows: $$G\times G/\Gamma \to G/\Gamma$$ $$(g,x\Gamma)\to gx\Gamma$$
One thinks of $G/\Gamma$ as the quotient space induced by the action on the right of $\Gamma$ on $G$, that is: $$\Gamma \times G\to G, \ \ (\gamma, g)\to g\gamma^{-1}$$
Now, let $(g^t)$ be a one-parameter subgroup of $G$, and consider the flow generated by the action on the left of $(g^t)$ on $G$, $i.e$: $$\R\times G\to G, \ \ (t,x)\to g^tx$$ 
Let $x\in G$ and define the curve $\gamma_x(t)$ on $G$ as follows: $\gamma_x(t)=g^tx$ for every $t\in \R$. Hence $\gamma_x(t)=R_x(g^t)= R_x(\gamma_1(t))$ where $1\in G$ is the identity element and $R_x$ is the right multiplication by $x$. It is clear now that: $$\gamma_x^{\prime}(0)=D_1R_x(\gamma_1^{\prime}(0))$$
Thus $\gamma_1^{\prime}(0)\in T_1G=\sf{Lie}(G)$ and the flow defined by the multiplication by $(g^t)$ is obtained by integrating the right invariant vector field $X(x)=D_1R_x(\gamma_1^{\prime}(0))$ for every $x\in G$. Conversely, if $X$ is a right invariant vector field on $G$, and $(g^t)$ is the one-parameter group for which $\left.\frac{\partial}{\partial t}  \right|_{t=0}g^t=\gamma_1^{\prime}(0)=X(1)$ where $\gamma_1(t)=g^t$, then the flow generated by $X$ is nothing but the left action of $(g^t)$ on $G$.\\\\
One may ask the following question: what are all the diffeomorphisms $\widetilde{\varphi}:G\to G$ preserving the set of all the right invariant vector fields on $G$? That is, if $X$ is a right invariant vector field, then the vector field $Y=\widetilde{\varphi}_*(X)$ is also right invariant?\\\\
Since any right invariant vector field $X$ corresponds to a left action on $G$ by a one-parameter subgroup $(g^t)$, then the question is in fact equivalent to finding all the diffeomorphisms $\widetilde{\varphi}$ that send a left action by a one-parameter subgroup to a left action by another one-parameter subgroup. More precisely, for every one-parameter subgroup $(g_1^t)$ there exists another one $(g_2^t)$ such that: $$\widetilde{\varphi}\circ L_{g_1^t}\circ \widetilde{\varphi}^{-1}=L_{g_2^t}$$
where $L_g$ denotes the left multiplication by $g$. In particular, $\widetilde{\varphi}$ preserves (by conjugacy) the set of left multiplications. Before answering the question, let us consider the following examples of such diffeomorphisms:
\begin{itemize}
\item \textit{Automorphisms:} If $\widetilde{\varphi}:G\to G$ is a Lie group automorphism, then: $$\widetilde{\varphi}\circ L_g\circ \widetilde{\varphi}^{-1}(x)=\widetilde{\varphi}(g\widetilde{\varphi}^{-1}(x))=\widetilde{\varphi}(g)\widetilde{\varphi}(\widetilde{\varphi}^{-1}(x))=L_{\widetilde{\varphi}(g)}(x)$$
\item \textit{Right translations:} Assume $\widetilde{\varphi}=R_h$, then: $$\widetilde{\varphi}\circ L_g\circ \widetilde{\varphi}^{-1}(x)=R_h(gxh^{-1})=gx=L_g(x)$$
\item \textit{Left translations:} If $\widetilde{\varphi}=L_h$, then: $$\widetilde{\varphi}\circ L_g\circ \widetilde{\varphi}^{-1}=L_h\circ L_g\circ L_{h}^{-1}=L_{hgh^{-1}}$$
\end{itemize}
\begin{proposition}
Let $\widetilde{\varphi}:G\to G$ be a diffeomorphism preserving (by conjugacy) the set of all the left translations. Then, $\widetilde{\varphi}=A\circ R_h$ where $A\in \sf{Aut}(G)$.
\end{proposition}
Proof: Consider the map $A=\widetilde{\varphi}\circ R_{h^{-1}}$ where $h^{-1}=\widetilde{\varphi}^{-1}(1)$. Then, it is enough to show that $A\in \sf{Aut}(G)$. We have that $A$ preserves the set of all the left multiplications because both $\widetilde{\varphi}$ and $R_{h^{-1}}$ do. Let $f:G\to G$ defined as follows: $A\circ L_g \circ A^{-1}=L_{f(g)}$. Then:$$ A\circ L_g\circ L_h\circ A^{-1}=(A\circ L_g\circ A^{-1})\circ (A\circ L_h\circ A^{-1})=L_{f(g)}\circ L_{f(h)}=L_{f(g)f(h)}$$
On the other hand we have: $$A\circ L_g\circ L_h\circ A^{-1}=A\circ L_{gh}\circ A^{-1}=L_{f(gh)}$$
Hence, $f(gh)=f(g)f(h)$, $i.e$ $f$ is a homomorphism. On the other side: $$A\circ L_g\circ A^{-1}(1)=L_{f(g)}(1)$$ 
This implies that $A(g)=f(g)$. So $A\in \sf{Aut}(G)$.\begin{flushright}
				$\square$
			\end{flushright}
\begin{remark}
We know that left multiplications preserve (by conjugacy) left multiplications. So, if $\widetilde{\varphi}=L_g$ then $\widetilde{\varphi}=A\circ R_g$ where $A\in \sf{Aut}(G)$. In fact, $A$ is the inner automorphism $L_gR_{g^{-1}}$.
\end{remark}
Let now $\mathcal{F}$ be a framing of $G$ obtained by right translating a basis $\mathcal{B}$ of $\sf{Lie}(G)=T_1G$. Hence, $\mathcal{F}$ defined by the $n$ right invariant vector fields determined by the $n$ vectors of $\mathcal{B}$ where $n=\dim(G)$. Let $A\in \sf{Aut}(G)$. Since $A$ preserves the set of all the right invariant vector fields, then $A$ is autonomous with respect to $\mathcal{F}$ with cocycle matrix $M_A$ equals the matrix of $D_1A:\sf{Lie}(G)\to \sf{Lie}(G)$ where $\sf{Lie}(G)$ is identified to $\R^n$ using the basis $\mathcal{B}$.\\\\
Consider now $g\in G$ and let $\sf{Ad}(g)$ denotes the inner automorphism $L_gR_{g^{-1}} \in \sf{Aut}(G)$. We know that $\sf{Ad}(g)$ is autonomous with respect to the right invariant framing $\mathcal{F}$. But we have that $R_{g^{-1}}$ preserves $\mathcal{F}$. This implies that $L_g:G\to G$ is autonomous with respect to $\mathcal{F}$ and with cocycle matrix $M_{L_g}=M_{\sf{Ad}(g)}$.\\\\
In fact, if $\widetilde{\varphi}:G\to G$ is a diffeomorphism which is autonomous with respect to a right invariant framing $\mathcal{F}$, then $\widetilde{\varphi}$ preserves the set of all the right invariant vector fields on $G$, which implies that $\widetilde{\varphi}$ preserves (by conjugacy), the set of all the left multiplications. Thus, using the previous proposition, we conclude that:
\begin{corollary}
Let $\widetilde{\varphi}:G\to G$ be an autonomous diffeomorphism with respect to a right invariant framing $\mathcal{F}$. Then $\widetilde{\varphi}=A\circ R_h$ where $A\in \sf{Aut}(G)$.
\end{corollary}
Let us investigate now the case of autonomous diffeomorphisms of the quotient space $G/\Gamma$. But before this, it is useful to remind some general facts. Consider a diffeomorphism $\varphi:G/\Gamma\to G/\Gamma$ and let $\varphi(\Gamma)=g\Gamma$ (that is, the point $g\Gamma\in G/\Gamma$ is the image of the point $\Gamma$ by the diffeomorphism $\varphi$). Define the new diffeomorphism $\psi=g^{-1}\varphi:G/\Gamma\to G/\Gamma$ (recall that $G$ acts on the left on $G/\Gamma$). Hence, $\psi$ fixes the point $\Gamma\in G/\Gamma$.\\\\ 
Consider now the action of $\psi$ on the space of all the curves $\gamma:[0,1]\to G/\Gamma$ such that $\gamma(0)=\Gamma\in G/\Gamma$. Let $\pi :G\to G/\Gamma$ be the natural projection $i.e$ the universal cover of $G/\Gamma$. Then for any curve $\gamma:[0,1]\to G/\Gamma$ such that $\gamma(0)=\Gamma$ there is a unique lifting $\widetilde{\gamma}:[0,1]\to G$ of $\gamma$ and satisfies $\widetilde{\gamma}(0)=1\in G$. Thus, in this way, $G$ is identified with the space of all the homotopy classes of the curves on $G/\Gamma$ starting from the point $\Gamma$. Let $\gamma$ be such a curve, then it corresponds to a unique point in $G$, that is, the point $\widetilde{\gamma}(1)$ where $\widetilde{\gamma}$ is the unique lifting of $\gamma$ satisfying $\widetilde{\gamma}(0)=1$. Define the map $\widetilde{\psi}:G\to G$ as follows: $$\widetilde{\psi}(x)=\widetilde{\psi\circ \gamma}(1)$$ where $\widetilde{\psi\circ \gamma}$ is the lifting of $\psi\circ \gamma$. The map $\widetilde{\psi}$ preserves $\Gamma$, and it is clear that $\widetilde{\psi}$ is bijective and covers the diffeomorphism $\psi$, $i.e$ $\pi\circ \widetilde{\psi}=\psi\circ \pi$. Hence, $\widetilde{\psi}$ is a diffeomorphism, and we obtain that $\varphi=g\psi$ is covered by the diffeomorphism $\widetilde{\varphi}=g\widetilde{\psi}:G\to G$.\\\\
Coming back to our situation, let us assume that $\mathcal{F}$ is a framing of $G/\Gamma$ induced by projecting a right invariant framing $\widetilde{\mathcal{F}}$ of $G$ using the natural map $\pi:G\to G/\Gamma$ (recall that this is well defined since the quotient space $G/\Gamma$ is obtained by the action of $\Gamma$ on the right). Let $\varphi:G/\Gamma \to G/\Gamma$ be an autonomous diffeomorphism with respect to the framing $\mathcal{F}$ and let $\widetilde{\varphi}=g\widetilde{\psi}: G\to G$ be the diffeomorphism that covers $\varphi$ as constructed above. Since $\pi \circ \widetilde{\varphi}=\varphi\circ \pi$ then $\widetilde{\varphi}$ is autonomous with respect to $\widetilde{\mathcal{F}}$ with cocycle matrix $M_{\widetilde{\varphi}}=M_{\varphi}$ (because the map $\pi$ transfers, locally, what happens in $G$ to what happens in $G/\Gamma$). This implies that $g^{-1}\widetilde{\varphi}=\widetilde{\psi}$ is autonomous with respect to $\widetilde{\mathcal{F}}$ such that $\widetilde{\psi}(1)=1$ and $\widetilde{\psi}(\Gamma)=\Gamma$. So, using the previous corollary, we get $\widetilde{\psi}\in \sf{Aut}(G)$. As a conclusion we have:
\begin{corollary}
Let $\mathcal{F}$ be a framing of $G/\Gamma$ induced by projecting a right invariant framing $\widetilde{\mathcal{F}}$ of $G$. Assume $\varphi:G/\Gamma\to G/\Gamma$ is an autonomous diffeomorphism with respect to $\mathcal{F}$. Then there exists $A\in \sf{Aut}(G)$ satisfying $A(\Gamma)=\Gamma$, and $g\in G$ such that $\varphi$ is covered by the affine map $L_g\circ A$.
\end{corollary}
\begin{remark}
Recall that the cocycle matrix $M_{L_g}=M_{\sf{Ad}(g)}$ and $M_A=D_1A:\sf{Lie}(G)\approx \R^n\to \sf{Lie}(G)$ with respect to $\mathcal{F}$. Hence, we have that $M_{\varphi}=M_{\sf{Ad}(g)}M_A$.
\end{remark}
\subsection{Integration of infinitesimal actions of Lie algebras} Let $G$ be a connected Lie group with a Lie algebra $\sf{Lie}(G)$ seen as the algebra of right invariant vector fields. Assume that $\psi:G\times M\to M$ is a smooth action of $G$ on $M$. For any right invariant vector field $X\in \sf{Lie}(G)$, there exists a unique one-parameter subgroup $(g^t)=\exp(tX)$ of $G$ whose tangent vector at $t=0$ is the vector $X(1)$. So, this defines an action of $\R$ (a flow) on $M$ as follows: $$\R\times M\to M $$ $$(t,x)\to \psi(\exp(tX(1)),x)$$ This flow is obtained by integrating the vector field: $$X_M(x)=\left.\frac{\partial}{\partial t}  \right|_{t=0} \psi(\exp(tX(1)),x)$$
That is, $X_M(x)\in T_x M$ is the tangent vector of the curve $\psi(\exp(tX(1)),x)$ at time $t=0$. This defines a map $$d\psi: \sf{Lie}(G)\to \Gamma^{\infty}(M), \ \ \ X\to X_M$$ where $\Gamma^{\infty}(M)$ is the space of all the smooth vector fields of $M$. The map $d\psi$ is a homomorphism of Lie algebras and the action $\psi$ is uniquely determined by $d\psi$, that is, if $\psi^{'}$ is another action of $G$ on $M$ such that $d\psi=d\psi^{'}$, then $\psi=\psi^{'}$. We call $d\psi$ the infinitesimal action.\\\\
It is natural to ask whether a homomorphism of Lie algebras $\sf{Lie}(G)\to \Gamma^{\infty}(M)$ corresponds to an infinitesimal action of the Lie group $G$?
\begin{theorem}[Palais 1957]\label{palais} (For a sketch of proof in the $C^1$ case see Section  \ref{appendix}) Let $G$ be a simply connected Lie group with Lie algebra $\sf{Lie}(G)$ of right invariant vector fields. Let $M$ be a manifold and assume that $\phi:\sf{Lie}(G)\to \Gamma^{\infty}(M)$ is a Lie algebra homomorphism. Assume furthermore that $\phi(X)\in \Gamma^{\infty}(M)$ is complete for every $X\in \sf{Lie}(G)$. Then $\phi$ is the infinitesimal action of a unique action of $G$ on $M$.

\end{theorem}

Assume $M$ is a compact connected manifold of dimension $n$ and $\mathcal{F}$ is a smooth framing of $M$ generated by the smooth vector fields $X_1, X_2, \ldots, X_n$. Assume, in addition, that: $$[X_i,X_j]=\sum_{k=1}^{n}a_k^{ij} X_k$$
for every $i, j\in \{1,2,\ldots,n\}$ where $a_1^{ij},a_2^{ij}, \ldots,a_n^{ij}$ are constant real numbers. Let $\mathcal{G}\subset \Gamma^{\infty}(M)$ be the vector sub-space generated by $X_1, X_2, \ldots, X_n$. Then $\mathcal{G}$ is a $n-$dimensional real Lie sub-algebra. Let $G$ be the associated $n-$dimensional real simply connected Lie group. Since the inclusion $i:\mathcal{G}=\sf{Lie}(G)\to \Gamma^{\infty}(M)$ is a Lie algebra homomorphism, then, according to the Palais's theorem, this corresponds to the infinitesimal action of a unique action of $G$ on $M$.\\\\
The condition that $X_1, X_2, \ldots, X_n$ define a basis of $T_x M$ at any $x\in M$ implies that the action of $G$ on $M$ is locally free. Thus each orbit of this action is open. This implies that $G$ acts transitively on $M$ which gives that $M=G/\Gamma$ where $\Gamma=\sf{Stab}(x)$ is the stabilizer of a chosen point $x\in M$. In addition, the vector fields $X_1, X_2, \ldots, X_n$ on $M=G/\Gamma$ correspond to the projection of right invariant vector fields of $G$. Combining with the previous results we obtain:
\begin{corollary}
Assume $M$ is a compact connected manifold of dimension $n$ and $\mathcal{F}$ is a smooth framing as indicated above. Let $\varphi: M\to M$ be an autonomous diffeomorphism with respect to $\mathcal{F}$. Then $M=G/\Gamma$ and $\varphi$ is covered by an affine map $L_g\circ A:G\to G$ with $A\in \sf{Aut}(G)$ satisfying $A(\Gamma)=\Gamma$ (where $G$, $\Gamma$, and the identification $M=G/\Gamma$ are as discussed above).
\end{corollary} 

\section{\textbf{Proof of Theorem \ref{2D} in the non-parabolic case}}
Assume that $M$ is a two dimensional compact surface and a framing $\mathcal{F}$ given by two vector fields $X$ and $Y$ (this implies that $M$ is diffeomorphic to the torus). Let $\varphi$ be an autonomous diffeomorphism with matrix: $$M_{\varphi}= \begin{pmatrix}
a &b \\ 
c &d 
\end{pmatrix}$$
We have seen (in \ref{properties}) that $\det(M_{\varphi})=\pm 1$. Let $[X,Y]$ be the Lie bracket of $X$ and $Y$, and $\varphi_*([X,Y])$ be its push forward by $\varphi$. Then we have:$$\varphi_*([X,Y])=[aX+cY,bX+dY]=(ad-bc)[X,Y]=\pm [X,Y]$$
In particular, if $\varphi$ preserves the orientation then $[X,Y]$ is $\varphi-$invariant, that is,  $\varphi_*([X,Y])=[X,Y]$.
\subsection{The hyperbolic case:} 
\begin{proposition}
Consider a framing $\mathcal{F}$ given by two vector fields $X$ and $Y$ and assume that $M_{\varphi}$ is hyperbolic. Then $[X,Y]=0$.
\end{proposition}
Proof: We have that $Z=[X,Y]$ is $\varphi-$invariant (as observed above). Let $x\in \T^2$ and suppose that $Z(x)\neq 0$. So, since $M_{\varphi}$ is hyperbolic, either $D\varphi^n(Z(x))$ or $D\varphi^{-n}(Z(x))$ goes to infinity. But this contradicts the fact that $Z$ is $\varphi-$invariant and bounded.\begin{flushright}
				$\square$
			\end{flushright}
\begin{corollary}
The diffeomorphism $\varphi$ is smoothly conjugate to an Anosov linear automorphism of $\T^2$.
\end{corollary}
Proof: Because in this case we have a transitive locally free action of $\R^2$ on the torus $\T^2$ generated by the flows of the vector fields $X$ and $Y$. Hence $\T^2=\R^2/\Gamma$ and $\mathcal{F}$ is obtained by projecting a parallel framing $\widetilde{\mathcal{F}}$ of $\R^2$. Furthermore, there is a diffeomorphism $\widetilde{\varphi}:\R^2\to\R^2$ covering $\varphi$ and has a constant derivative (with respect to the framing $\widetilde{\mathcal{F}}$) and equals $M_{\varphi}$. Hence $\widetilde{\varphi}$ is affine. But since $\widetilde{\varphi}$ is hyperbolic, then it must have at least one fixed point. This implies that it is conjugate by a translation to its linear part.\begin{flushright}
				$\square$
			\end{flushright}

\subsection{The elliptic case}
Suppose that the cocycle matrix $M_{\varphi}$ of $\varphi$ with respect to $\mathcal{F}$ is elliptic. Then $M_{\varphi}$ belongs to the isometry group of a Euclidean inner product on $\R^2$, denoted by $q_0$. Since each tangent space $T_x \T^2$ is identified to $\R^2$ via the basis $\mathcal{F}(x)$, then, we can define the smooth Riemannian metric $g$ on $\T^2$ by inducing $q_0$ on each tangent plane via the identification. In this case, $\varphi$ is an isometry of the smooth Riemannian metric $g$.
\begin{proposition}
The metric $g$ is conformally equivalent to a flat metric.
\end{proposition}
Proof: This is exactly the content of the uniformization theorem (see for instance \cite{Hubb} for more details). More precisely, $g$ defines an almost complex structure which is integrable (since the dimension is two). The uniformization theorem implies that this complex torus is the quotient of $\C$ by a lattice of translations. Therefore $g$ is conformally flat.
\begin{flushright}
				$\square$
			\end{flushright}
\begin{corollary}
The map $\varphi$ is an automorphism of a complex torus, hence, it is affine.
\end{corollary}
Proof: Let $f:\T^2\to \R_{>0}$ be smooth map such that the metric $fg$ is flat. We have that $\varphi$ is conformal with respect to the metric $fg$, because it is an isometry of $g$. The flat metric $fg$ defines a unique complex structure of $\T^2$ for which $\varphi$ is biholomorphic. Automorphisms of complex tori are known to be affine (see \cite{Hubb} for more details).\begin{flushright}
				$\square$
			\end{flushright}
\begin{remark}
In this case we have $\T^2=\C/\Gamma$ and $\varphi$ is covered by an automorphism $\widetilde{\varphi}:\C\to\C$ of the form $z\to az+b$ such that $a\Gamma=\Gamma$. This implies that $a\in \S^1$ (and rational). Hence $\varphi$ is, in fact, an isometry of the flat metric.
\end{remark}

\section{\textbf{Proof of Theorem \ref{2D} in the parabolic case}}
Up to changing the framing, we can suppose that $\mathcal{F}$ is a smooth framing of $\T^2$ for which the diffeomorphism $\varphi$ is autonomous with cocycle matrix: $$ M_{\varphi}=\begin{pmatrix}
1 &1 \\ 
0 &1 
\end{pmatrix}$$
The framing $\mathcal{F}$ is given by two independent vector fields $X$ and $Y$. Thus both $X$ and $[X,Y]$ are $\varphi-$invariant.
\begin{proposition}
We have $[X,Y]=\beta X$ such that $\beta$ is a smooth $\varphi-$invariant function defined on $\T^2$.
 \end{proposition}
 Proof: If a vector field $Z$ is $\varphi-$invariant and for some $a\in \T^2$ we have $Z(a)\neq \lambda X(a)$ then the norm of the iterates of $Z(a)$ under $\varphi$ (in the associated Riemannian metric) goes to infinity. This is impossible because $Z$ is $\varphi-$invariant. The fact that $\beta$ is $\varphi-$invariant comes from the fact that $[X,Y]$ is $\varphi-$invariant.\begin{flushright}
				$\square$
			\end{flushright}
 \begin{lemma}
 Suppose that $\beta$ is constant non null. Then there is a locally free transitive action of the affine Lie group $\sf{Aff}(\R)$ on $\T^2$.
 \end{lemma}
 
 Proof: We can suppose that $\beta=1$ (up to changing the framing), that is, $[X,Y]=X$. This defines a finite Lie subalgebra of the algebra of all the smooth vector fields generated by $X$ and $Y$ isomorphic to the Lie algebra of $\sf{Aff}(\R)$. Hence, applying the Palais' theorem, we conclude that this defines an action of $\Aff(\R)$ on $\T^2$ which is transitive and locally free (because $X$ and $Y$ define a framing).\begin{flushright}
				$\square$
			\end{flushright}
 \begin{corollary}
 If the function $\beta$ is constant then it vanishes and $\varphi$ is smoothly conjugate to an affine parabolic automorphism of the torus.
 \end{corollary}
 Proof: The fact that "if $\beta$ is constant then it must vanish" comes from the fact that there is no locally free transitive action of the group $\Aff(\R)$ on $\T^2$ because $\Aff(\R)$ doesn't contain any discrete subgroup $\Gamma$ for which $\Aff(\R)/\Gamma$ is compact (since it is not unimodular). Thus we are in the flat case, which we already investigated.\begin{flushright}
				$\square$
			\end{flushright}
 \begin{lemma}\label{4.3}
 We have that $d\beta(X)=0$, that is, the line field generated by $X$ is contained in $\ker(d\beta)$.
 \end{lemma}
 Proof: Let $z\in \T^2$ be a recurrent point for $\varphi$. Suppose that $X(z)\notin \ker(d_z\beta)$ that is $X(z)$ is transversal to the line $\ker(d_z\beta)$. In this case, we have that the line $\ker(d_z\beta)$ converges under the iterates of $\varphi$ to be tangent to $X$ because $\varphi$ is autonomous with respect to $X$ and $Y$ with cocycle matrix $ M_{\varphi}=\begin{pmatrix}
1 &1 \\ 
0 &1 
\end{pmatrix}$. But on the other side $\varphi$ preserves $\ker(d\beta)$ because $\beta$ is $\varphi-$invariant. So, since $z$ is recurrent, we conclude that $\ker(d_z\beta)$ is tangent to $X$. Which contradicts our assumption. Hence $X(z)\in \ker(d_z\beta)$. But almost every point is recurrent (because $\varphi$ preserves a finite Riemannian volume), this implies $d\beta(X)=0$.\begin{flushright}
				$\square$
			\end{flushright}

Suppose now that $\beta$ is not constant, then for a generic point $a$ in the image of $\beta$, the level $\beta^{-1}(a)$ (which is $\varphi-$ and $X-$invariant) is a compact one dimensional manifold ($i.e$ union of circles) tangent to the vector field $X$.\\\\
The equality $[X,Y]=\beta X$ implies that the flow generated by $Y$ sends $X-$orbits to $X-$orbits (not necessarily preserving parametrizations). More accurately, the flow $\phi^t_Y$ preserves the $X-$foliation.\\\\
We want to show that for a generic $\beta-$level, its saturation by the $Y-$flow covers the whole $\T^2$.
 \begin{proposition} Let $F$ be a codimension one foliation on a compact manifold $N$ with at least one compact leaf. assume there is a flow $\phi^t$ transversal to $F$ and sends leaf to leaf. Then $F$ corresponds to a fibration $N\to \S^1$.
\end{proposition}
Proof: Let $L_0$ be a compact leaf. Hence $L_0$ is a compact embedded submanifold of $M$. Define the map $i:L_0\times \R\to M$ by $i(x,t)=\phi^t(x)$. Since $\phi^t$ is transversal to $F$ and sends leaf to leaf then the map $i$ is an immersion and for every $t\in \R$ there is $\epsilon >0$ such that the restriction of $i$ on $L_0\times]t-\epsilon, t+\epsilon[$ is an embedding (because $L_0$ is compact). The map $i$ cannot be injective, in fact if $L_x$ is any leaf, define $$O_x=\{ \phi^t(y) \mid t\in \R, y\in L_x\}$$
then $O_x$ is open (because $\phi^t$ is transversal to $F$) and either $O_{x_1}\cap O_{x_2}=\emptyset$ or $O_{x_1}= O_{x_2}$ (because $\phi^t$ sends leaf to leaf). So we have a partition of $M$ into open sets. But $M$ is connected, and this implies that we have only one open set. Hence the map $i$ is surjective. If furthermore it is injective then it is a diffeomorphism which is impossible. Let $t_0\in \R$ be defined as follows:$$t_0=\inf \{t>0 \mid \phi^t(L_0)=L_0\}$$
So $t_0>0$ because $i$ is a local diffeomorphism and $t_0\neq \infty$ otherwise $i$ is injective. Thus $M\approx L_0\times [0,t_0]$ with $(x,0)\sim (\phi^{t_0}(x), t_0)$.\begin{flushright}
				$\square$
			\end{flushright}
\begin{remark}
Note that the assumption that the foliation has at least one compact leaf is crucial. It is needed to assure that the image of this compact leaf $L_0$ by the transversal flow, in a small neighbourhood of the time $0$, defines a foliation by disjoint leaves in a neighbourhood of $L_0$. A counterexample without the compactness assumption is the foliation of the torus by dense parallel lines with a transversal translation flow. 
\end{remark}
\begin{corollary}
In this case $\varphi$ preserves a fibration $F$ by circles tangent to the vector field $X$.
\end{corollary}
Proof: We have that $\varphi$ preserves the foliation $F$ defined by $X$. Since $\beta$ is not constant, a generic level $\beta^{-1}(a)$ is a finite union of circles tangent to $X$. Hence, applying the previous proposition we conclude that $F$ is a fibration by circles all tangent to $X$ and preserved by $\varphi$.\begin{flushright}
				$\square$
			\end{flushright}
 \begin{lemma}\label{conj}
 Let $f:\S^1\to \S^1$ be an orientation preserving diffeomorphism that preserves a framing of $\S^1$, then $f$ is smoothly conjugate to the rotation $R_{\alpha}$ (where $\alpha$ is the rotation number of $f$).
 
 \end{lemma}
 Proof: This happens because $f$ is an isometry of the associated metric. In fact let $x\in\S^1$ and let $\gamma:[0,1]\to \S^1$ be the simple arc joining $x$ and $f(x)$ and moving in the direction of the framing. Let $l(\gamma)$ be the length of $\gamma$ in the associated metric and $d$ be the total length of $\S^1$. Put $\alpha=2\pi(l(\gamma)/d)$. Then clearly the isometry $f$ is conjugate to $R_{\alpha}$.\begin{flushright}
				$\square$
			\end{flushright}
Let $X$ be a smooth non-vanishing vector field on $\S^1\times\S^1$ such that $X(z_1,z_2)$ is tangent to the circle $\{z_1\}\times \S^1$ for every $(z_1,z_2)\in \S^1\times \S^1$. Denote by $(\S^1\times \S^1, X)$ to mean the space $\S^1\times \S^1$ with the vector field $X$ on it, and let $$\varphi:(\S^1\times \S^1,X)\to (\S^1\times \S^1,X)$$ be a diffeomorphism preserving $X$ and preserves each factor $\{z\}\times \S^1=S_z$. Hence, for each $z\in\S^1$, the restriction of $\varphi$ on $S_z$ is smoothly conjugate to the associated rotation, denoted $R_{\alpha(z)}$. In this case we have the following:

\begin{lemma}
There exists a global diffeomorphism $\phi: (\S^1\times \S^1,X)\to \S^1\times \S^1$ preserving each factor $\{z\}\times \S^1=S_z$ such that the restriction of $\phi\circ \varphi\circ \phi^{-1}$ on each $S_z$ equals $R_{\alpha(z)}$.
\end{lemma}
 Proof: Let $(z_1,z_2)\in \S^1\times \S^1$ and let $l(z_1,z_2)$ be the length of the unique curve starting from $(z_1,1)$ and ending at $(z_1,z_2)$ and moving in the direction of the restriction of the vector field $X$ on $S_z$ (where the length is computed with respect to the metric given by $X$ on $S_z$). Let $d(z_1)$ be the total length of $S_{z_1}$ with respect to this metric. Define $f_{z_1}: \{z_1\}\times \S^1=\S^1\to \S^1$ by: $$f_{z_1}(z_2)=e^{i(l(z_2)/d(z_1))2\pi}$$
 Put now $\phi(z_1,z_2)=(z_1, f_{z_1}(z_2))$. From the construction we observe that $f_{z_1}$ sends the restriction of $X$ on $\{z_1\}\times \S^1$ to the usual unitary framing of $\S^1$. Hence $\phi$ sends the vector field $X$ to the constant unitary vector field on $\S^1\times \S^1$ tangent to the fibers $\{z\}\times \S^1$ and in the positive orientation on each fiber (the fibers are identified with $\S^1$). This implies that $\phi$ satisfies the desired condition.\begin{flushright}
				$\square$
			\end{flushright}
Now let $\varphi:\S^1\times \S^1\to \S^1\times \S^1$ be a diffeomorphism, preserving each factor $\{z\}\times \S^1$, and $X$ be a non-vanishing vector field on $\S^1\times \S^1$ all as defined above. Let $\alpha:\S^1\to \S^1$ be the map satisfying that the restriction of $\varphi$ on each $\{z\}\times\S^1$ is smoothly conjugate to the associated rotation $R_{\alpha(z)}$ (by Lemma \ref{conj}). From the previous lemma, we observe that $\alpha$ is smooth. We have:

\begin{lemma}
If the map $\alpha:\S^1\to \S^1$ is a local diffeomorphism (hence a covering) then $\varphi$ is smoothly conjugate to the linear automorphism on $\R^2/ \Z^2$ generated by the matrix $ \begin{pmatrix}
1 &k \\ 
0 &1 
\end{pmatrix}$ where $k$ is the degree of this covering.
\end{lemma}
Proof: Let $z_1, z_2, \ldots, z_k$ be the points of the fiber $\alpha^{-1}(1)$ ordered in the positive direction. Let $\gamma:[0,1]\to \S^1$ be the curve defined by $\gamma(t)= z_1 e^{(2\pi t)i}$. Hence $\gamma$ is a simple closed curve such that $\gamma(0)=\gamma(1)=z_1$. Put $\pi:\R\to \S^1$ to be the universal covering given by $\pi(x)=e^{(2\pi x)i}$ and let $\widetilde{\alpha\circ \gamma}:[0,1]\to \R$ be the lifting of the curve $\alpha\circ \gamma$ such that $\widetilde{\alpha\circ \gamma}(0)=0$. Thus $\widetilde{\alpha\circ \gamma}(1)=k$. In addition we have $\widetilde{\alpha\circ \gamma}(t_i)=i-1$ for each $t_i$ satisfying $\gamma(t_i)=z_i$, \ $1\leq i \leq k$. Define the map $f:\S^1\to \R/\Z$ as follows: $$f(z)=\frac{\widetilde{\alpha\circ \gamma}(\gamma^{-1}(z))}{k} \mod (1)$$
note that this map is well-defined because $\gamma^{-1}(z_1)=\{0,1\}$ has two choices but they are the same $\mod (1)$. Consider the map (as defined in the proof of the previous lemma) $l:\S^1\times \S^1\to\R_{\geq 0}$ such that $l(z_1,z_2)$ is the length of the unique curve starting from $(z_1,1)$ and ending at $(z_1,z_2)$ and moving in the direction of the restriction of the vector field $X$ on $S_z$. Similarly, let $d:\S^1\to \R_{>0}$ such that $d(z)$ is the total length of $\{z\}\times \S^1$. The maps $l$ and $d$ are smooth. Define finally the map $\phi: \S^1\times \S^1 \to \R/\Z\times\R/\Z$ as follows: $$\phi(z_1,z_2)=(f(z_1), \frac{l(z_1,z_2)}{d(z_1)}) \mod (1)$$
Hence $\phi\circ \varphi\circ \phi^{-1}$ is equal to $ \begin{pmatrix}
1 &k \\ 
0 &1 
\end{pmatrix}$ because one sees that the restriction of the diffeomorphism $\phi\circ \varphi\circ \phi^{-1}$ on each segment $\{t\}\times[0,1]$ for $t\in[0,1]$ is exactly the vertical translation by the value $kt \mod (1)$.\begin{flushright}
				$\square$
			\end{flushright}
\begin{remark}
If the map $\alpha$ is supposed to be only a covering map (without being a local diffeomorphism) then $\varphi$ is still conjugate to a linear parabolic automorphism of $\R^2/\Z^2$. But, in this case, the conjugacy is just $C^0$ $i.e$ the conjugating map is just a homeomorphism.
\end{remark}
Let us now come back to our situation. Suppose $\varphi:\T^2\to\T^2$ be an autonomous diffeomorphism with respect to a smooth framing $\mathcal{F}$, given by the vector fields $X$ and $Y$, and having a parabolic cocycle matrix $M_{\varphi}=\begin{pmatrix}
1 &1 \\ 
0 &1 
\end{pmatrix}$. We have seen that $\varphi$ preserves a fibration $F$ by circles tangent to the vector field $X$. The fibration $F$ is equivalent (by a diffeomorphism) to the trivial fibration of $\S^1\times\S^1$ given by the vertical circles $\{z\}\times \S^1$, $z\in \S^1$. Hence the quotient $\T^2/F$ is diffeomorphic to $\S^1$ and $\varphi$ induces a diffeomorphism $\widehat{\varphi}:\T^2/F\to \T^2/F$. The flow $\phi_Y^t$ generated by $Y$ preserves $F$, $i.e$ sends circle to circle. This implies that $Y$ projects to a smooth vector field $\widehat{Y}$ on $\T^2/F$. Furthermore, $\widehat{\varphi}$ preserves $\widehat{Y}$ because we have: $$\widehat{\varphi}_*(\widehat{Y})=\pi_*(\varphi_*(Y))=\pi_*(X+Y)$$
Where $\pi:\T^2\to \T^2/F$ is the quotient map. So, since $X$ is tangent to $F$, we have $\pi_*(X)=0$ which implies $\pi_*(X+Y)=\pi_*(Y)=\widehat{Y}$. Hence $\widehat{\varphi}_*(\widehat{Y})=\widehat{Y}$ and, consequently, $\widehat{\varphi}$ is smoothly conjugate to a rotation.\\\\
Since the function $\beta$ is not constant (recall that $[X,Y]=\beta X$ and $\beta$ is $\varphi-$invariant) then $\varphi$ preserves a generic fiber $\beta^{-1}(a)$ which is a finite union of circles tangent to $F$. This implies that $\widehat{\varphi}$ is smoothly conjugate to a periodic rotation of period $n\geq 1$. Thus the diffeomorphism $\varphi^n:\T^2\to \T^2$ preserves each circle of the fibration $F$ and has cocycle matrix $M_{\varphi^n}=\begin{pmatrix}
1 &n \\ 
0 &1 
\end{pmatrix}$ with respect to the framing $\mathcal{F}$. In this case we have:
\begin{lemma}
Let $\varphi^n$ be the diffeomorphism as above and let $R:\T^2/F\to \S^1$ be the map associating the rotation number of the restriction of $\varphi^n$ on each circle of the fibration $F$. Then $R$ is a local diffeomorphism.
\end{lemma}
Proof: Let $\alpha:\T^2\to \R$ be a positive non vanishing smooth function which is $\varphi-$invariant, that is, $\alpha$ is constant along each leaf of the fibration $F$. Define $X^{\prime}=\alpha X$ and $Y^{\prime}=\alpha Y$. Then it is clear that $\varphi^n$ is still autonomous with the framing given by $X^{\prime}$ and $Y^{\prime}$ and having the same parabolic cocycle matrix because $\varphi^n$ preserves each fiber of the fibration $F$. In addition, we have that: $$[X^{\prime},Y^{\prime}]=[\alpha X,\alpha Y]=\alpha^2\beta X-\alpha Y(\alpha)X=\alpha(\alpha\beta-Y(\alpha))X=(\alpha\beta-Y(\alpha))X^{\prime}$$
where $[X,Y]=\beta X$ and we have that $(\alpha\beta-Y(\alpha))$ is constant along each fiber just as we discussed before in Lemma \ref{4.3}. So $[X^{\prime},Y^{\prime}]=0$ if and only if $\alpha\beta-Y(\alpha)=0$. Let now $p\in \T^2$ and $\gamma:]-\epsilon,\epsilon[\to \T^2$ be an embedded smooth small curve such that $\gamma(0)=p$ and tangent to $Y$. Then $\beta$ restricted to $\gamma$ is a well defined function and the differential equation $\alpha\beta-Y(\alpha)=0$ on $\gamma$ has a local solution around $p$, that is, there exists a unique $\alpha_{x_0}:\gamma(]-\epsilon^{\prime},\epsilon^{\prime}[)\to \R$ with $\epsilon^{\prime}<\epsilon$ satisfying $\alpha_{x_0}(p)=x_0$ for $x_0\in\R_{>0}$. Then $\alpha_{x_0}$ is positive around $p$. Consequently, we can construct a global function $\alpha:\T^2\to \R_{>0}$ which is constant along each fiber of $F$ and coincides with $\alpha_{x_0}$ around $p$. Since $(\alpha\beta-Y(\alpha))$ is constant along fibers, then it vanishes on a cylinder $C$ containing the fiber of $p$. So, the restriction of the framing given by $X^{\prime}$ and $Y^{\prime}$ on the cylinder $C$ is flat, which implies that $C$ has a translation structure where the flow given by $X^{\prime}$ is by translations. consider the segment $I_p$ in $C$ which is around $p$ and tangent to $Y^{\prime}$, then its image under the $X^{\prime}-$flow comes back to itself identically since each point is periodic and the $X^{\prime}-$flow on $C$ sends a line to a line parallel to it (with respect to the translation structure). Thus $C$ is universally covered by a strip band $B=]a,b[\times\R$ with a translation acting on it. This implies that the diffeomorphism $\varphi^n$ restricted to $C$ can be lifted to diffeomorphism of $B$ preserving each level $\{x\}\times\R$ and having the same cocycle matrix  with respect to the canonical framing. Consequently, the function on $C$ associating rotation number of $\varphi^n$ restricted to each fiber of $F$ in $C$ is a local diffeomorphism because this is reduced to the linear case (locally).

\begin{flushright}
				$\square$
			\end{flushright}
As observed by Pierre Mounoud, we can, with another argument, verify that the
differential equation in the proof of the previous lemma has a global solution on $\mathbb{S}^1$. This simplifies the situation and reduces the problem to the flat case.

 \begin{corollary}
The diffeomorphism $\varphi^n$ is smoothly conjugate to the parabolic linear automorphism of $\R^2/\Z^2$ given by the matrix $ A_k=\begin{pmatrix}
1 &k \\ 
0 &1 
\end{pmatrix}$ where $k$ is the degree of the covering $R:\T^2/F\to \S^1$.

\end{corollary}\begin{flushright}
				$\square$
			\end{flushright}
This completes the proof of Theorem \ref{2D}. Observe that all we needed in the proof is the fact that $\beta$ is $C^2$ (to deduce the existence of regular levels of $\beta$ using Sard). Therefore, the same conclusion holds assuming only $C^3$ regularity of the framing.
\section{\textbf{Proof of Theorem \ref{3D}}}
\subsection{$C^{\infty}-$Regularity:} For the sake of clarity and ease of reading we present in this section a proof assuming smooth regularity. Then, in the next section \ref{C1}, we deal with the $C^1$ case. So, assume that everything is smooth, particularly, the vector fields fields $X_s, X_c$ and $X_u$
\begin{lemma}
Assume that $\lambda_c\neq \pm 1$, then $[X_i,X_j]=0$ for $i,j\in \{s,c,u\}$.
\end{lemma}
Proof: Let $i,j\in \{s,c,u\}$. We have that $\varphi_*([X_i,X_j])=[\lambda_iX_i,\lambda_jX_j])=\lambda_i\lambda_j[X_i,X_j]$. Suppose that $[X_i,X_j]\neq 0$ and let $x\in M$ such that $[X_i,X_j](x)\neq 0$. Then the vector field $[X_i,X_j]$ does not vanish on an open neighbourhood $V$ of $x$. Put $U$ to be the open subset obtained as the saturation of $V$ (the union of all the $\varphi-$iterates of $V$). Then $[X_i,X_j]$ does not vanish on $U$ and the claim is that $[X_i,X_j]$ (restricted to $U$) is proportional to one of the vector fields $X_s, X_c$, or $X_u$. Indeed, let $l$ be the line field on $U$ generated by $[X_i,X_j]$, then $l$ is $\varphi-$invariant and for a generic recurrent point $y\in V$ the iterates of the line $l(y)$ come back sufficiently close to $l(y)$. This forces $l(y)$ to be tangent to one of the directions $X_s, X_c$, or $X_u$. Since we can choose $V$ small enough, we have that $[X_i,X_j]$ is proportional to one of the vector fields $X_s, X_c$, or $X_u$ and, so, the same holds on $U$. Thus $[X_i,X_j]=\beta X_k$ for some function $\beta$ and $k\in \{s,c,u\}$. Furthermore, on $U$ we have: $$\varphi_*([X_i,X_j](x))=\lambda_i\lambda_j[X_i,X_j](\varphi(x))=\lambda_i\lambda_j\beta(\varphi(x))X_k(\varphi(x))$$
On the other hand:
$$\varphi_*(\beta(x) X_k(x))=\lambda_k\beta(x) X_k(\varphi(x))$$  
So, $\lambda_i\lambda_j\beta(\varphi(x))=\lambda_k\beta(x)$, or, in other words:

$$\beta(\varphi(x))=\frac{\lambda_k}{\lambda_i\lambda_j}\beta(x)$$

So $\beta$ must be identically zero on $U$ since $\left|\frac{\lambda_k}{\lambda_i\lambda_j}\right|\neq1$ and almost every point recurrent.\begin{flushright}
				$\square$
			\end{flushright}
\begin{corollary}
If $\lambda_c\neq \pm 1$ then $\varphi$ is smoothly conjugate to a linear Anosov automorphism of the three dimensional torus $\T^3=\R^3/\Z^3$.
\end{corollary}
Proof: The framing given by $X_s, X_c$, and $X_u$ is flat and the cocycle matrix is hyperbolic. This implies the desired conclusion.\begin{flushright}
				$\square$
			\end{flushright}
			\begin{remark}\label{Remark 5.3}
Note that, in this case, the same results still hold if we assume only that the vector fields $X_s, X_c$, and $X_u$ are just $C^1$. The only difference is that the conjugacy in this situation is $C^1$.
			
			\end{remark}
It remains to treat the case $\lambda_c=\pm 1$. But we can restrict ourselves to the case $\lambda_c=1$ and $\lambda_s,\lambda_u>0$ by considering $\varphi^2=\varphi\circ \varphi$. So, from now on, we suppose $\lambda_c=1$, $\lambda_s,\lambda_u>0$.
\begin{lemma}
We have $[X_c,X_s]=\beta X_s$ and $[X_c,X_u]=\alpha X_u$ where $\beta$ and $\alpha$ are smooth $\varphi-$invariant functions on $M$.
\end{lemma}
Proof: We have that $\varphi_*([X_c,X_s])=[X_c,\lambda_s X_s]=\lambda_s[X_c,X_s]$. Hence, the vector field $[X_c,X_s]$ is contracting under the action of $\varphi$ which implies that $[X_c,X_s]=\beta X_s$ where $\beta$ is a smooth function on $M$. We have that $D_x\varphi(\beta(x)X_s(x))=\lambda_s\beta(x)X_s(\varphi(x))$. On the other hand, since $\varphi_*([X_c,X_s])=\lambda_s[X_c,X_s]$, we have that: $$D_x\varphi(\beta(x)X_s(x))=\lambda_s\beta(\varphi(x))X_s(\varphi(x))$$
So $\beta(x)=\beta(\varphi(x))$. The same holds for $X_c$ and $X_u$ by considering $\varphi^{-1}$.\begin{flushright}
				$\square$
			\end{flushright}
   Throughout the remaining of this section, we will treat the case where $\beta$ is not constant. We prove in this case (Corollary \ref{Lor}) that up to finite power and cover, $\varphi$ is smoothly conjugate to either $A\times Id$ on $\T^2\times \S^1=\T^3$, or to an affine automorphism $L_g\circ A:G/\Gamma\to G/\Gamma$ where  $G=\sf{Lor}(1,1)=\SO(1,1)\ltimes \R^2$. The case where $\beta$ is constant will be treated in the next section \ref{algebraic}.
\begin{proposition}\label{tangent}
Suppose that $\beta$ is not constant. Then for a generic regular point $a$ in the image of $\beta$ we have that the level surface $\beta^{-1}(a)$ is tangent to the plane distribution given by the vector fields $X_s$ and $X_u$.

\end{proposition}
\begin{lemma}
If a continuous function $f:M\to Y$ is $\varphi-$invariant then $f$ is invariant by the flows $\phi_s^t$ and $\phi_u^t$ generated by $X_s$ and $X_u$ respectively.
\end{lemma}
Proof: Let $x\in M$ and $\gamma:[0,1]\to M$ be a $C^1$ curve tangent to the vector field $X_s$ with $\gamma(0)=x$ and $\gamma(1)=y$. Let $l(\gamma)$ be the length of this curve with respect to the Riemannian metric given by the framing. Since $\gamma^{'}(t)=\alpha(t)X_s(\gamma(t))$ we obtain that $(\varphi\circ \gamma)^{'}(t)=D_{\gamma(t)}(\gamma^{'}(t))=\lambda_s\gamma^{'}(t)$. This implies that $l(\varphi\circ\gamma)=\lambda_s l(\gamma)$. Hence the curve $\gamma$ is getting smaller under the iterates of $\varphi$ which implies that $x$ and $y$ are getting closer. But $\beta$ is $\varphi-$invariant, thus, $\beta(\varphi^n(x))=\beta(x)$ and $\beta(\varphi^n(y))=\beta(y)$. Since $\varphi^n(x)$ and $\varphi^n(y)$ keep getting closer and closer, we obtain $\beta(x)=\beta(y)$ (using the compactness of $M$). Consequently, $\beta$ is $\phi_s^t-$invariant and the same argument holds for the unstable flow $\phi_u^t$ by considering $\varphi^{-1}$ instead. This proves the lemma and the proposition since the level surface $\beta^{-1}(x)$ is, hence, $\phi_s^t$ and $\phi_u^t$ invariant. This proves Proposition. \ref{tangent}\begin{flushright}
				$\square$
			\end{flushright}

We know that $\varphi$ preserves the level surface $\beta^{-1}(a)$. Suppose, furthermore, that $\varphi$ preserves each connected component of $\beta^{-1}(a)$ (otherwise consider a finite power of $\varphi$). Thus, $\varphi$ preserves a two dimensional torus $\T^2\subset M$ tangent to the plane distribution given by $X_s$ and $X_u$ (by Proposition \ref{tangent}). Consequently, $\varphi$ acts as a linear Anosov diffeomorphism on $\T^2$. In this case we have:
\begin{proposition}\label{product}
The image of this torus $\T^2$ by the central flow $\phi_c^t$ defines a fibration of $M$ over $\S^1$ where each fiber is diffeomorphic to $\T^2$ and preserved by $\varphi$. Furthermore, $\phi_c^t$ sends fiber to fiber by conjugating the dynamics of $\varphi$ on each one.
\end{proposition}\label{susp}
Proof: Consider the saturation of this torus $\T^2$ by the central flow. This is clearly an open set $U$ since $\T^2$ is transversal to $\phi_c^t$. The set $U$ is also invariant by the flows $\phi_s^t$ and $\phi_u^t$ because $\phi_c^t$ preserves the orbits of $\phi_s^t$ and $\phi_u^t$ respectively and $\T^2$ is tangent to the stable and the unstable directions. But the subsets of $M$ saturated by $\phi_c^t$, $\phi_s^t$, and $\phi_u^t$ define a partition of $M$ by open subsets. So we have $U=M$ and the torus $\T^2$ must come back since $M$ is compact.  Hence the image of the torus $\T^2$ by the central flow $\phi_c^t$ defines a fibration of $M$ over $\S^1$ where each fiber is diffeomorphic to $\T^2$. Thus, we have a projection $p:M\to \S^1$ where $\S^1$ denotes the space of fibers on which $\varphi$ induces an autonomous diffeomorphism with respect to the projection of the central vector field. So, this induced diffeomorphism is conjugate to a rotation, but, it has a fixed point (corresponding to the invariant torus). Consequently, $\varphi$ preserves each fiber and, since $\varphi$ commutes with the central flow, then $\phi_c^t$ sends fiber to fiber by conjugating the dynamics of $\varphi$ on each one.

\begin{flushright}
				$\square$
			\end{flushright}
Let $A:\T^2\to \T^2$ be an Anosov linear automorphism and consider the diffeomorphism $\widetilde{A}:\T^2\times [0,1]\to \T^2\times [0,1]$ defined by $\widetilde{A}(x,t)=(A(x),t)$. Let $f:\T^2\to \T^2$ be a diffeomorphism that commutes with $A$ and consider $M$ to be the three-manifold defined by identifying $(x,0)$ to $(f(x),1)$ in the product space $\T^2\times [0,1]$. Then $\widetilde{A}$ induces a diffeomorphism $A_f$ on $M$. 
\begin{corollary} \label{5.6}
The diffeomorphism $\varphi$ (up to finite power) is smoothly conjugate to some $A_f$ for $A:\T^2\to \T^2$ an Anosov linear automorphism and $f:\T^2\to \T^2$ a diffeomorphism that commutes with $A$.
\end{corollary}
Proof: This follows immediately from Proposition \ref{susp}.
\begin{flushright}
				$\square$
			\end{flushright}

\begin{corollary}\label{Lor}
    Up to finite power and cover, $\varphi$ is smoothly conjugate to either $A\times Id$ on $\T^2\times \S^1=\T^3$, or to a partially hyperbolic affine automorphism $L_g\circ A:G/\Gamma\to G/\Gamma$ where  $G=\sf{Lor}(1,1)=\SO(1,1)\ltimes \R^2$.
\end{corollary}
Proof: Let $A$, $f$, and $A_f$ be as in the statement of Corollary \ref{5.6}. Let $f_t\subset\SL(2,\R)$ be the one parameter group of hyperbolic matrices such that $f_1=f$.\\
Suppose now that there is $k\in \Z$ such that $f=A^k$. If $k=0$ then $f=Id$ and $\varphi$ is smoothly conjugate to $A\times Id$ on $\T^2\times \S^1=\T^3$. Assume $k\neq 0$ and let $G=\R\ltimes\R^2$ where $t\in \R$ acts via $f_t$. Clearly $G=\sf{{Lor}}(1,1)$. We have $A\in (f_t)$  and we see it as an element $\widetilde{A}$ of $G$ via the inclusion $(f_t)\subset G$. We have in this case that $\Ad(\widetilde{A})$ acts on $G$ as a partially hyperbolic automorphism acting on the normal factor $\R^2\subset G$ via $A$ itself. Let $\Gamma\subset G$ be the cocompact lattice isomorphic to $\Z\ltimes\Z^2$ where $\Z\subset (f_t)$ generated by $f$ and acts on $\R^2$ preserving $\Z^2$. Then $\Ad(\widetilde{A})$ preserves $\Gamma$ and induces a partially hyperbolic diffeomorphism on $G/\Gamma$ which is by construction smoothly conjugate to $A_f$.\\
The point now, is that the centralizer of an Anosov automorphism of the two dimensional torus consists, up to finite index, of the iterates of the automorphism (see for instance \cite{PY} and \cite{Pl}).
In higher dimension $n$, we can have  actions of $\Z^{n-1}$
centralizing an Anosov automorphism. So the result is not true, yet it is true, generically! as shown in \cite{PY}. Returning to our situation, this shows that up to finite cover $f$ is an iterate of $A$ and, therefore, the previous discussion applies to this finite cover.

\begin{flushright}
				$\square$
			\end{flushright}

\subsection{Algebraic classification}\label{algebraic} In this section we deal with the case where the functions $\alpha$ and $\beta$ are constant.
\begin{proposition}
We have that $[X_s,X_u]=\delta X_c$ where $\delta$ is a constant real number.
\end{proposition}
Proof: We have $\varphi([X_s,X_u])=[X_s,X_u]$. So the vector field $[X_s,X_u]$ must be tangent to the central direction (otherwise it cannot be preserved). Hence $[X_s,X_u]=\delta X_c$ and $\delta$ is $\varphi-$invariant. Suppose $\delta$ is not constant and let $a$ in the image of $\delta$ be a regular value. Thus all the points in a neighbourhood of $a$ are regular values and their level sets define a foliation of an open subset $O=\delta^{-1}(]a-\epsilon, a+\epsilon[)$ of $M$. This implies that the plane field generated by $X_s$ and $X_u$ is integrable on $O$ which implies that $[X_s,X_u]$ is tangent to this plane field. But in the same time it is tangent to the central direction. Consequently $\delta$ equals zero on $O$ which contradicts the fact that $a$ is a regular value. Thus $\delta$ is constant.

\begin{flushright}
				$\square$
			\end{flushright}
We have seen in the previous subsection that $[X_c,X_s]=\beta X_s$ and $[X_c,X_u]=\alpha X_u$ and we understood the diffeomorphism $\varphi$ in the case where $\beta$ is not constant. The same holds if $\alpha$ is not constant. So it remains to classify in the case where $[X_c,X_s]=a_{cs} X_s$, $[X_c,X_u]=a_{cu} X_u$ and $[X_s,X_u]=a_{su}X_c$ such that $a_{cs}$, $a_{cu}$ and $a_{su}$ are constant real numbers.\\\\
In this case, assuming $\alpha$ and $\beta$ constant, the vector fields $X_s$, $X_c$, and $X_u$ generate a three dimensional real Lie subalgebra of the Lie algebra of all smooth vector fields. Put $\mathcal{G}=<X_s, X_c, X_u>$ to be the subalgebra spanned by these vector fields and let $G$ be the (unique) associated connected simply connected Lie group (that is, $\sf{Lie}(G)=T_1G\approx \mathcal{G}$). Palais' theorem implies that $G$ acts on $M$ transitively and locally freely (because $X_s$, $X_c$, and $X_u$ define a basis at each point of the compact manifold $M$) and that $M$ is identified with $G/\Gamma$ where $\Gamma$ is a cocompact lattice ($i.e$ $\Gamma$ is a discrete subgroup of $G$ with $G/\Gamma$ compact). In addition, the vector fields $X_s$, $X_c$, and $X_u$ are the projection of right invariant vector fields $\widetilde{X_s}$, $\widetilde{X_c}$, and $\widetilde{X_u}$ of $G$ on the quotient $G/\Gamma=M$. Thus our diffeomorphism $\varphi: M=G/\Gamma\to G/\Gamma$ is covered by an affine automorphism $\widetilde{\varphi}:G\to G$, that is, $\widetilde{\varphi}=L_g\circ A$ where $L_g$ is a left translation on $G$ and $A\in \sf{Aut}(G)$ satisfying $A(\Gamma)=\Gamma$. The cocycle matrix $M_{\varphi}$ of $\varphi$ with respect to the framing $\mathcal{F}$ generated by $X_s$, $X_c$, and $X_u$ is equal to the cocycle matrix $M_{\widetilde{\varphi}}$ of $\widetilde{\varphi}=L_g\circ A$ with respect to the framing $\widetilde{\mathcal{F}}$ generated by $\widetilde{X_s}$, $\widetilde{X_c}$, and $\widetilde{X_u}$. And we have: $M_{\widetilde{\varphi}}=M_{\sf{Ad}(g)}M_A$ (where $M_{\sf{Ad}(g)}$ and $M_A$ are the cocycle matrices of $\sf{Ad}(g)$ and $A$ respectively).\\\\
So, our situation is the following: we have a three-dimensional connected simply connected Lie group $G$ and $\Gamma\subset G$ is a discrete subgroup with $G/\Gamma$ is compact. Furthermore, we have $A\in \sf{Aut}(G)$ and $g\in G$ such that $A(\Gamma)=\Gamma$ and the matrix $M_{\sf{Ad}(g)}M_A$ is partially hyperbolic with respect to a basis $(X_s, X_c, X_u)$ of the algebra $\sf{Lie}(G)$ of right invariant vector fields. In this case we have a classification.\\\\
It is known that if $G$ is a simply connected Lie group of dimension three having a cocompact lattice $\Gamma\subset G$, then $G$ must be one of the following six groups (for more details see, for example, the paper \cite{RV}):
\begin{itemize}
\item The Abelian group $\R^3=\R \oplus \R \oplus \R$.
\item The nilpotent group $\sf{Heis}_3$ (the Heisenberg group) of $3\times 3$ upper triangular real matrices of the form $\begin{pmatrix}
1 &x &z \\ 
0 &1 &y \\
0 &0 &1 
\end{pmatrix}$.
\item The solvable group $\widetilde{\sf{Euc}}(\R^2)$: the universal cover of the Euclidean group of orientation preserving isometries of $\R^2$ ($i.e$ $\sf{SO}(2)\ltimes \R^2$).
\item The solvable group $\sf{Lor}(1,1)=\SO(1,1)\ltimes \R^2$ (also called $\sf{SOL}$): the identity component of the isometry group (the Lorentzian group) of the Minkowski space $(\R^2, q)$ where $q$ is the quadratic form $q(x,y)=xy$. It is isomorphic to the semi-direct product $\R\ltimes \R^2$ where $\R$ acts on $\R^2$ via the representation $t\to \begin{pmatrix}
e^t &0 \\ 
0 &e^{-t} 
\end{pmatrix}$.
\item The simple compact group $\sf{SU}(2)\approx \S^3$.
\item The simple non-compact group $\widetilde{\sf{SL}}(2,\R)$: The universal cover of the special linear group $\sf{SL}(2,\R)$. This is the richest case.
\end{itemize}
\begin{lemma} Let $G$ be either $\sf{SU}(2)$ or $\widetilde{\sf{Euc}}(\R^2)$.
Let $\widetilde{\varphi}=L_g\circ A:G\to G$ be an affine automorphism. Then the cocycle matrix of $\widetilde{\varphi}$ cannot be partially hyperbolic with respect to a right invariant vector field.

\end{lemma}
Proof: If $G=\sf{SU}(2)$ then the Lie algebra of $G$ is isomorphic to $\R^3$ endowed with the bilinear map given by the usual cross product. The automorphisms of this algebra are exactly the linear Euclidean isometries. Hence $\sf{Aut}(G)$ is compact. Therefore the affine group $\sf{Aff}(G)=\{L_g\circ A \mid A\in \sf{Aut}(G), g\in G\}$ is compact. Let $\mathcal{F}$ be a right invariant framing and consider the cocycle representation: $$\Aff(G)\to \sf{Lie}(G)=\GL(3,\R)$$ $$L_g\circ A\to M_{\sf{Ad}(g)}M_A$$
Thus, the image of this representation is compact. So, it does not contain any partially hyperbolic matrix.\\\\
Let now $\widetilde{\sf{Euc}}(\R^2)=\R\ltimes \R^2$ where $\R$ acts on $\R^2$ via the representation $\rho :t\in\R\to \exp(tA)$. where $A$ is the matrix $\begin{pmatrix}
0 &-1 \\ 
1 &0 
\end{pmatrix}$.\\ 
Claim: all one-parameter subgroups $\R\subset \widetilde{\sf{Euc}}(\R^2)$ which are transversal to $\R^2$ (the subgroup of translations) are conjugate to each other. To see this we observe that such a one-parameter subgroup acts on $\R^2$ by rotations around some fixed point, hence, the conjugacy between two one-parameter subgroups transversal to $\R^2$ can be done by the translation between the two fixed points.\\
The Lie algebra $\mathcal{G}$ of $\widetilde{\sf{Euc}}(\R^2)$ is generated by $T,X,Y$ where $T$ corresponds to $\R$ and $X,Y$ correspond to $\R^2$ with bracket relations:$$ [T,X]=AX, \ \ [T,Y]=AY$$
Let $\varphi:\mathcal{G}\to \mathcal{G}$ be an automorphism such that $\varphi(T)=T$. Then:$$\varphi(AX)=\varphi([T,X])=[\varphi(T),\varphi(X)]=[T,\varphi(X)]=A(\varphi(X))$$
Hence $\varphi\circ A=A\circ \varphi$. Let now $\varphi:\mathcal{G}\to \mathcal{G}$ be partially hyperbolic. We can assume that $\varphi(T)=\alpha T$ (because all one-parameter subgroups transversal to $\R^2$ are conjugate). Suppose $\alpha\neq 1$, then:$$\varphi(AX)=\varphi([T,X])=[\varphi(T),\varphi(X)]=[\alpha T,\varphi(X)]=\alpha A(\varphi(X))$$ 
which is impossible because $\det (A)\neq 0$. Thus, $T$ is the central direction and $\R^2=E^s\oplus E^u$. Since $\varphi\circ A=A\circ \varphi$ we have (up to conjugacy) that $A=\begin{pmatrix}
a &0 \\ 
0 &b 
\end{pmatrix}$ with $a+b=0$. This implies that $\mathcal{G}$ is the Lie algebra of $\sf{Lor}(1,1)$, contradiction. So $\widetilde{\sf{Euc}}(\R^2)$ does not contain any partially hyperbolic automorphism. 

\begin{flushright}
				$\square$
			\end{flushright}
\label{5.2} So, it remains to classify in terms of the four remaining cases: $\R^3, \ \sf{Heis}_3, \ \sf{Lor}(1,1)$, and $\widetilde{\sf{SL}}(2,\R)$.\\
\begin{itemize}
\item \textbf{The Abelian group $\R^3$:} In this case any discrete cocompact subgroup $\Gamma\subset \R^3$ is equivalent to $\Z^3=\Z \oplus \Z \oplus\Z$. The automorphisms of $\R^3$ preserving the lattice $\Z^3$ are exactly the elements of $\sf{GL}(3,\Z)$. An affine automorphism $A+v$ with $A(\Z^3)=\Z^3$ has cocycle matrix equal to $A$ (with respect to the canonical basis). Hence $A+v$ induces a map on $\T^3=\R^3/\Z^3$ which is partially hyperbolic if and only if $A\in \sf{GL}(3,\Z)$ is partially hyperbolic. Suppose $A$ is such a map, then we have the splitting of $\R^3$ into the three directions $E_s$, $E_c$, and $E_u$ (the stable, central, and unstable directions). Suppose $E_c$ is associated to the eigenvalue $\lambda_c=\pm 1$ (otherwise $A$ is Anosov). In this case $E_c$ always meets the lattice $\Z^3$ in a point other than the origin. Hence the foliation of $\R^3$ by the parallel lines to $E_c$ projects to a fibration by circles on $\T^3$. Also, the foliation of $\R^3$ by the planes parallel to $E_s\oplus E_u$  projects to a foliation by $2-$dimensional tori of $\T^3$ since the matrix $A$ acting on $\Q^3$ admits a rational splitting into a rational line and a rational plane. This shows that, up to finite cover, $A$ is conjugate to a direct product $T\times Id$ on $\T^2\times \S^1$ where $T$ is an Anosov automorphism on $\T^2$.
\item \textbf{The Heisenberg group $\sf{Heis}_3$:}
Define, for $k\in \N$, the discrete subgroup $\Gamma_k\subset \sf{Heis}_3$: $$\Gamma_k=\left\lbrace  \begin{pmatrix}
1 &m &l/k \\ 
0 &1 &n \\
0 &0 &1 
\end{pmatrix} \mid m, n ,l \in \Z \right\rbrace$$
Then $\Gamma_k$ is a (uniform) lattice and it is known that every lattice of $\sf{Heis}_3$ is equivalent to exactly one $\Gamma_k$ (see for example \cite{BL-RL}). The quotient space $\sf{Heis}_3/\Gamma_k$ is diffeomorphic to a circle extension of the torus $\T^2$ where the circle correspond to the central direction in $\sf{Heis}_3$.\\\\
Let $C\approx \R$ be the center of $\sf{Heis}_3$. Then we have an exact sequence: $$1\to C\to \sf{Heis}_3\to \R^2\to 1$$
Let now $\sf{Aut}(\sf{Heis}_3)$ be the automorphism group and $\sf{Inn}(\sf{Heis}_3)$ be the subgroup of inner automorphisms. Hence, $\sf{Inn}(\sf{Heis}_3)\approx \sf{Heis}_3/C\approx \R^2$ and we obtain a natural map $\sf{Aut}(\sf{Heis}_3)\to \sf{GL}(2,\R)$ since $C$ is preserved by automorphisms. This defines a sequence: $$1\to \sf{Inn}(\sf{Heis}_3)\approx \R^2 \to \sf{Aut}(\sf{Heis}_3)\to \sf{GL}(2,\R)\to 1$$
This sequence is, in fact, exact and it splits, $i.e$, $\sf{Aut}(\sf{Heis}_3)\approx\sf{GL}(2,\R)\ltimes\R^2$.\\\\
It is also known (see \cite{BL-RL}) that we have a similar decomposition of the automorphism group of a (uniform) lattice: $$\sf{Aut}(\Gamma_k)\approx  \sf{GL}(2,\Z)\ltimes\Z^2$$ In particular, the group $\sf{Aut}(\Gamma_k)$ contains partially hyperbolic automorphisms (with central direction tangent to the center $C$).\\

\item \textbf{The Lorentzian group $\sf{Lor}(1,1)=\SO(1,1)\ltimes \R^2$:}
We have that $\sf{Lor}(1,1)=\R\ltimes \R^2$ where $\R$ acts on $\R^2$ via the representation $t\to\begin{pmatrix}
e^t &0 \\ 
0 &e^{-t} 
\end{pmatrix}$. Let $\Gamma\subset \sf{Lor}(1,1)$ be a uniform lattice. Then $\Gamma_0=\Gamma\cap \R^2$ is a uniform lattice of $\R^2$. And also the projection $p:\Gamma \to \sf{Lor}(1,1)/\R^2\approx \R$ is a lattice ($i.e$ discrete). Let $\Gamma_1=p^{-1}(p(\Gamma))$. We have that $\Gamma_0\subset \Gamma$ is normal and $\Gamma_1\ltimes \Gamma_0=\Gamma$ ($\Gamma_1$ acts by conjugacy on $\Gamma_0$).\\\\
Let $(h^t)$ be the one-parameter subgroup of $\sf{Lor}(1,1)=\R\ltimes \R^2$ consisting of diagonal matrices and let $h^{\alpha}\in (h^t)$ such that $h^{\alpha}$ is conjugate to an element of $\sf{SL}(2,\Z)$. Then $h^{\alpha}$ preserves some lattice $\Gamma_0$ of $\R^2$ and $\Gamma=\{h^{n\alpha}\}\ltimes \Gamma_0\approx\Z\ltimes \Z^2$ is a uniform lattice of $\sf{Lor}(1,1)$. Left translations on $\sf{Lor}(1,1)$ by elements of $(h^t)$ are partially hyperbolic and they induce partially hyperbolic autonomous diffeomorphisms on $\sf{Lor}(1,1)/\Gamma$.\\

\item \textbf{The group $\sf{\widetilde{SL}}(2,\R)$:}
Let $\pi:\widetilde{\sf{SL}}(2,\R)\to \sf{PSL}(2,\R)$ be the covering projection. Then $\pi_1(\PSL(2,\R))=\ker(\pi)\approx\Z$ is a central discrete subgroup. Let $\widetilde{\Gamma}\subset \widetilde{\sf{SL}}(2,\R)$ be a uniform lattice and $\Gamma=\pi(\widetilde{\Gamma})$. Our claim is that $\Gamma$ is a uniform lattice in $\PSL(2,\R)$ (this seems to be a folkloric fact but not easy to localise in the literature. So, let us give a hint).\\\\
We have a natural map $\rho:\widetilde{\Gamma}\setminus \widetilde{\sf{SL}}(2,\R)\to \Gamma\setminus \PSL(2,\R)$ and all the point is to show that $\Gamma$ is 
discrete (because in this case $\rho$ is a local diffeomorphism and $\widetilde{\Gamma}\setminus \widetilde{\sf{SL}}(2,\R)$ is compact which implies that $\Gamma\setminus \PSL(2,\R)$ is compact and $\rho$ is a covering map).\\\\
The first observation is that $\Gamma$ is discrete if and only if $\widetilde{\Gamma}\cap \ker(\pi)$ has finite index in $\ker(\pi)$, that is, there is $N\in \N$ such that $\sigma^{N}\in \widetilde{\Gamma}$ where $\sigma$ is a generator of $\ker(\pi)$. Put $X=\widetilde{\Gamma}\setminus \widetilde{\sf{SL}}(2,\R)$ then $\sigma$ acts on $X$ by right multiplication. Endow $\widetilde{\sf{SL}}(2,\R)$ with a left invariant Riemannian metric $m$. Then $m$ is well defined on $X$ and the right action of $\sigma$ (which is also a left multiplication because $\sigma$ is central) preserves $m$, that is, $\sigma\in \sf{Isom}(X,m)$. Since $X$ is compact then $\sf{Isom}(X,m)$ is compact. If $\sf{Isom}(X,m)$ is finite then for some power $N$ we have $\sigma^N=Id$ and this implies that $\sigma^N\in \widetilde{\Gamma}$. Suppose now that the closure of $\{\sigma^n, \ n\in \Z\}$ contains a one-parameter subgroup, then there is a corresponding one-parameter subgroup $(h^t)$ in $\sf{Isom}(\widetilde{\sf{SL}}(2,\R),m)$. One can choose $m$ such that $\sf{Isom}(\widetilde{\sf{SL}}(2,\R),m)=\widetilde{\sf{SL}}(2,\R)$ acting on itself by left translations. Hence, the one-parameter group $(h^t)$ is a subgroup of left multiplications which implies that $\widetilde{\Gamma}$ is centralized by all elements of $(h^t)$. In particular, $\widetilde{\Gamma}\subset \sf{centralizer}((h^t))$. But, the observation is that $\sf{centralizer}((h^t))=(h^t).Z$ where $Z=\ker(\pi)$ and $\pi: \widetilde{\sf{SL}}(2,\R)\to \PSL(2,\R)$ is the covering projection). We have that: $$(h^t).Z\setminus \widetilde{\sf{SL}}(2,\R)=(h^t)\setminus \PSL(2,\R)$$ which is non-compact and this contradicts the fact that $\widetilde{\Gamma}\subset \sf{centralizer}((h^t))$. We conclude that $\Gamma=\pi(\widetilde{\sf{SL}}(2,\R))$ is a uniform lattice of $\PSL(2,\R)$.\\\\
Let now $\Gamma\subset \PSL(2,\R)$ be a uniform lattice and $\varphi\in \Aut(\PSL(2,\R))$ preserving $\Gamma$. One fact is that $\sf{Inn}(\PSL(2,\R))\subset \Aut(\PSL(2,\R))$ has finite index, so, up to a finite power, we may assume that $\varphi\in \sf{Inn}(\PSL(2,\R))$. But an inner automorphism $L_g R_{g^{-1}}$ and the left translation $L_g$ have the same cocycle matrix with respect to a right invariant framing (which is the matrix of $\sf{Ad}(g)$) and it is partially hyperbolic if and only if $g$ is a hyperbolic element of $\PSL(2,\R)$. In particular, $\Gamma$ contains many hyperbolic elements and the associated inner automorphism preserves $\Gamma$ and induces a partially hyperbolic diffeomorphism on $\PSL(2,\R)/\Gamma$.\\\\
An other example is the left action of the one-parameter subgroup of diagonal matrices of $\PSL(2,\R)$ on $\PSL(2,\R)/\Gamma$ by partially hyperbolic left translations (this is no other than the geodesic flow of $\Gamma\setminus \H$).

\end{itemize}

\subsection{The $C^1-$regularity}\label{C1}

Assume now that the vector fields $X_s$, $X_c$, and $X_u$ are only $C^1$. We have, as before, that $[X_c,X_s]=\beta X_s$ and $[X_c,X_u]=\alpha X_u$ but in this case $\beta$ and $\alpha$ are only continuous $\varphi-$invariant functions on $M$. The case where $\left|\lambda_c\right|\neq 1$ is noted in Remark \ref{Remark 5.3}. So, in the remaining part, we suppose $\lambda_c=1$.
\subsubsection{\textbf{Non-constant structure coefficients:}}
Let $\mathcal{A}$ denote the set of all the closed subsets of $M$ that are $\varphi-$invariant and invariant under the action of the flows $\phi_s^t$ and $\phi_u^t$. Clearly if $a$ is a point in the image of $\varphi-$invariant function $\beta$, then the set $\beta^{-1}(a)$ belongs to $\mathcal{A}$ (for the same argument as in the $C^{\infty}$ regularity). In addition, the central flow $\phi_c^t$ preserves $\mathcal{A}$ because if $F\in \mathcal{A}$ then $\varphi(\phi_c^t(F))=\phi_c^t(\varphi(F))=\phi_c^t(F)$ and $\phi_c^t$ preserves the stable and the unstable foliations (because of the equalities $[X_c,X_s]=\beta X_s$ and $[X_c,X_u]=\alpha X_u$). Recall that $\phi_c^t$ and $\varphi$ commute since $\varphi_*(X_c)=X_c$.\\\\
Suppose now that $\beta$ is not constant and let $x\in M$. Let $D_x$ be the minimal subset containing $x$ and belongs to $\mathcal{A}$, that is, the intersection of all the elements of $\mathcal{A}$ containing $x$. It is clear that $D_x$ is a proper subset of $M$ since $\beta^{-1}(\beta(x))$ is in $\mathcal{A}$, containing $x$, and it is a proper subset of $M$ (because $\beta$ is not constant). The aim is to prove that $D_x$ is a compact $C^1$ embedded surface.\\\\
For this, let $\Lambda_x=\{t\in \R \mid \phi_c^t(D_x)=D_x\}$. Clearly $\Lambda_x$ is a closed subgroup of $\R$. Suppose that $\Lambda_x=\R$, then $D_x$ is invariant under the flows $\phi_s^t, \phi_c^t$, and $\phi_u^t$. This implies that $D_x=M$ but this is impossible since $D_x$ is a proper subset of $M$. Hence $\Lambda_x=d_x\Z$ where $d_x$ is a positive real number. Note that if for some $t\in\R$ we have $\phi_c^t(D_x)\cap D_x\neq \emptyset$ then $\phi_c^t(D_x)=D_x$ because $\phi_c^t(D_x)\cap D_x$ belongs to $\mathcal{A}$ since both $\phi_c^t(D_x)$ and $D_x$ do, and $D_x$ is minimal.
\begin{lemma}
For every $x\in M$, the subset $D_x$ is a compact embedded $C^1$ surface tangent to the plane field generated by $X_s$ and $X_u$. 
\end{lemma}
Proof: Let $y\in D_x$ be any point and let $I_{\epsilon}=]-\epsilon,\epsilon[$. Let $\gamma_x^s(I_{\epsilon})$ be the image of the curve $\gamma$ tangent to $X_s$ and satisfies $\gamma(0)=x$, that is, $\gamma_x^s(I_{\epsilon})=\{\phi_s^t(x) \mid -\epsilon<t<\epsilon\}$. Let now: $$ S_x^{\epsilon}=\{\phi_u^t(\gamma_x^s(I_{\epsilon})) \mid -\epsilon < t <\epsilon \}$$
Then for $\epsilon$ small enough, $S_x^{\epsilon}$ is a small piece of a $C^1$ surface tangent to $X_s$ and $X_u$ contained in $D_x$. The central flow $\phi_c^t$ is transversal to $S_x^{\epsilon}$ and for a small time $\alpha$ the sets $\phi_c^t(S_x^{\epsilon})$ for $-\alpha <t<\alpha$ define a $C^1$ foliation of an open subset $B_x$ around $x$. But for $\alpha<d_x$ the subsets $\phi_c^t(S_x^{\epsilon})$ are all disjoint and do not intersect $D_x$ except $S_x^{\epsilon}$. This shows $S_x^{\epsilon}$ is the only piece of $D_x$ inside $B_x$ which proves the lemma.\begin{flushright}
				$\square$
			\end{flushright}
\begin{corollary}
The plane field generated by $X_s$ and $X_u$ is integrable and defines a $C^1$ foliation by compact leaves.
\end{corollary}
Proof: For any $x\in M$ the subset $D_x$ is a closed embedded $\varphi-$invariant surface tangent to $X_s$ and $X_u$.\begin{flushright}
				$\square$
			\end{flushright}
Consider, as before, $D_x$ to be a closed minimal $\varphi-$invariant and preserved by the flows $\phi_s^t$ and $\phi_u^t$. We proved that $D_x$ is a closed compact $C^1$ embedded surface. It is clear that each connected component of $D_x$ is $C^1$ diffeomorphic to the torus $\T^2$, because its tangent bundle is necessarily trivial. The diffeomorphism $\varphi$ acts by permuting this components in a cyclic transitive way, otherwise $D_x$ wouldn't be minimal. Let $C$ be such a component, then a finite power of $\varphi$ preserves $C$ (and each other component). The central flow $\phi_c^t$ preserves the foliation generated by the vector fields $X_s$ and $X_u$. Since $\phi_c^t$ is transversal to each leaf, then the orbit of every leaf is an open subset of $M$. But $M$ is connected, hence, The image of $C$ by central flow covers $M$. Following the same idea in the proof of proposition (\ref{product}) we conclude that $M$ is $C^1$ diffeomorphic to a suspension of $\T^2$, in other words, the foliation generated by $X_s$ and $X_u$ corresponds to a $C^1$ fibration $M\to \S^1$ with $\T^2$ as a typical fiber.  More precisely, we have exactly as in Corollary \ref{5.6} and Corollary \ref{Lor}:
\begin{corollary}
A finite power of the diffeomorphism $\varphi$ is $C^1$ conjugate to some $A_f$ for $A:\T^2\to \T^2$ an Anosov linear automorphism and $f:\T^2\to \T^2$ a diffeomorphism that commutes with $A$. Therefore, the same conclusions of Corollary \ref{Lor} apply in this $C^1$ setting.
\end{corollary}
Proof: Let $C$ be as indicated above (so $C$ is minimal). Then $C$ is $C^1$ diffeomorphic to $\T^2$ and preserved by a finite power $\varphi^k$ of $\varphi$. The action of $\varphi^k$ on $C$ is clearly hyperbolic and, hence, $C^1$ conjugate to a linear Anosov automorphism of $\T^2$. Furthermore, the image of $C$ by the central flow defines a $C^1$ foliation of $M$ that corresponds to a fibration over $\S^1$ as follows: let $t_0>0$ be the smallest time for which $\phi_c^{t_0}(C)=C$ (note that $t_0$ exists otherwise we would have an embedding of $C\times \R$. Which is surjective because the orbit of each leaf is open and since $M$ is connected we must have only one orbit). Then the action of the central flow on the space of leaves is periodic, with period $t_0$, and this implies that this corresponds to a fibration and the space of leaves is $C^1$ diffeomorphic to $\S^1$ on which the diffeomorphism $\varphi^k$ acts by preserving the framing induced by projecting $X_c$, hence it acts isometrically with respect to the corresponding metric. Since $\varphi^k$ preserves $C$ then it acts on the space of leaves trivially (preserving each leaf) because this action is isometric (since it preserves the projection of $X_c$) with a fixed point. Thus $\varphi^k$ preserves each leaf of the fibration. In addition, the central flow $\phi_c^t$ commutes with $\varphi^k$, hence $\phi_c^t$ sends leaf to leaf by conjugating the induced dynamics on each one. In particular, $\phi_c^{t_0}$ commutes with the restriction of $\varphi^k$ on $C$ (where $t_0$ is, as mentioned before, the smallest time for which $\phi_{t_0}^c(C)=C$). This proves the result.\begin{flushright}
				$\square$
			\end{flushright}
\subsubsection{\textbf{Case of constant structure coefficients:}} We have the equality: $$\varphi_*([X_s,X_u])=[\lambda_sX_s,\lambda_uX_u]=[X_s,X_u]$$
Hence $[X_s,X_u]$ is $\varphi-$invariant and, for the same reasons as before, this implies that $[X_s,X_u]$ is tangent to the central direction $i.e$ $[X_s,X_u]=\delta X_c$ where $\delta$ is a $\varphi-$invariant continuous function on $M$. If $\delta$ is not constant, then we are in the previous situation which we understood. So, it remains to treat the case where the functions $\alpha$, $\beta$, and $\delta$ are constant. Put $\alpha=a_{cu}$, $\beta=a_{cs}$, and $\delta=a_{su}$. Thus, we have: $[X_s,X_u]=a_{su}X_c$, $[X_c,X_s]=a_{cs}X_s$, and $[X_c,X_u]=a_{cu}X_u$ where $a_{su}, a_{cs}, a_{cu}\in \R$. Then we have:
\begin{lemma}
The vector fields $X_s$, $X_c$, and $X_u$ generate a three-dimensional Lie algebra, that is, the Jacobi identity is satisfied.
\end{lemma}
\begin{remark}
Observe that we cannot consider derivation of brackets since the vector fields are only $C^1$.
\end{remark}
Proof: Note that it is enough to verify the identity just for the triplet $X_s$, $X_c$, and $X_u$. We have the following equality: $$[[X_s,X_c],X_u]+[[X_c,X_u],X_s]+[[X_u,X_s],X_c]=-a_{su}(a_{cs}+a_{cu})X_c$$
So, in order to have the Jacobi identity, we must have either $a_{su}=0$ or $a_{cs}+a_{cu}=0$. But the central flow is autonomous and preserves the stable and the unstable directions and for some time $t_0$ it has determinant equal to $e^{(a_{cs}+a_{cu})}$ which must equal $1$. Hence $a_{cs}+a_{cu}=0$. \begin{flushright}
				$\square$
			\end{flushright}
In this case, we have that the $C^1$ vector fields $X_s$, $X_c$, and $X_u$ generate an abstract three dimensional real Lie algebra denoted $\mathcal{G}$. Let $G$ be the associated connected simply connected Lie group. So, (see \ref{appendix}) there is a unique $C^1$ action of $G$ on the compact manifold $M$ such that the embedding of $\sf{Lie}(G)=\mathcal{G}$ into the space of $C^1$ vector fields on $M$ sending $\mathcal{G}$ onto $<X_s,X_c,X_u>$ is seen as the corresponding infinitesimal action. More precisely, $G$ acts (a $C^1$ action) on $M$ transitively and locally freely such that $X_s$, $X_c$, and $X_u$ are the projections of right invariant vector fields $\widetilde{X_s}$, $\widetilde{X_c}$, and $\widetilde{X_u}$ of $G$ on the quotient $M=G/\Gamma$ where $\Gamma$ is the stabilizer of a chosen point. The autonomous diffeomorphism $\varphi:M=G/\Gamma\to M$ is induced by an affine map $L_g\circ A:G\to G$ with $A\in \sf{Aut}(G)$ satisfying $A(\Gamma)=\Gamma$. As a conclusion we have:
\begin{corollary}
The diffeomorphism $\varphi:M=G/\Gamma\to M$ is covered by an affine automorphism $L_g\circ A:G\to G$ such that $A\in \sf{Aut}(G)$ satisfying $A(\Gamma)=\Gamma$. Thus, we have the same classification as in the $C^{\infty}$ regularity, the only difference is that the conjugacy is just $C^1$.
\end{corollary}
\section{\textbf{Counterexamples to Theorem \ref{3D} in the $C^0-$regularity case}}\label{C0} We have seen that in the $C^1$ regularity case we have an algebraic classification of partially hyperbolic autonomous dynamical systems. The situation when assuming only the $C^0$ regularity of the framing is different as previously discussed in Section \ref{1.1}. The present section is devoted to give more details by constructing autonomous partially hyperbolic diffeomorphisms (called skew-product of Anosov diffeomorphism) which are, generically, non-algebraic (this implies, according to Theorem \ref{3D}, that the framing is not $C^1$).\\\\
Let $\varphi:M\to M$ be a diffeomorphism and let $\rho:M\to G$ be a $C^k$ map where $G$ is a compact connected Lie group. Define the map: $$\varphi_{\rho}:M\times G\to M\times G, \ \ \varphi_{\rho}(x, g)=(x,\rho(x)g)$$
So, $\varphi_{\rho}$ sends the fiber $\{x\}\times G$ to the fiber $\{\varphi(x)\}\times G$ by applying the left translation $L_{\rho(x)}$ on the second factor $i.e$ on $G$. In particular, we have $\pi\circ \varphi_{\rho}=\varphi\circ \pi$ where $\pi:M\times G\to M$ is the projection on $M$. The diffeomorphism $\varphi_{\rho}$ is called a $G-$extension (or a skew-product) of $\varphi$. If $G=\S^1$ then it is simply called a "circle extension" of $\varphi$.\\\\
Suppose that $\mathcal{F}_G$ is a left invariant framing on $G$. Then we can define a family of independent global vector fields on $M\times G$ by inducing $\mathcal{F}_G$ on each fiber $\{x\}\times G$. Furthermore, these vector fields are preserved by $\varphi_{\rho}$ because it acts by left translations on the fibers.\\\\
Suppose now that $\varphi:M\to M$ is Anosov and $\varphi_{\rho}:M\times G\to M\times G$ as before. Suppose in addition that there exists a horizontal $\varphi_{\rho}-$invariant distribution $\mathcal{H}$ of $T(M\times G)$. That is, $\mathcal{H}$ is a sub-bundle of $T(M\times G)$ such that $\mathcal{H}(x) \oplus T_xG_x=T_x(M\times G)$ and $D_x\varphi_{\rho}(\mathcal{H}(x))=\mathcal{H}(\varphi_{\rho}(x))$ for every $x\in M\times G$ where $G_x$ is the fiber $\{\pi(x)\}\times G$. In this case the diffeomorphism $\varphi_{\rho}$ is partially hyperbolic such that the $\varphi-$invariant splitting $E_s\oplus E_u=TM$ is lifted to $\varphi_{\rho}-$invariant splitting $\widetilde{E}_s\oplus \widetilde{E}_u=\mathcal{H}$ and the Whitney sum: $$\widetilde{E}_s\oplus \widetilde{E}_c\oplus \widetilde{E}_u=T(M\times G)$$ satisfies the partial hyperbolicity conditions (where $\widetilde{E}_c=\ker(d\pi)$ is the vertical distribution tangent to the fibers $G_x$, $x\in G$).
\begin{fact}
A circle extension $\varphi_{\rho}:M\times \S^1\to M\times \S^1$ of an Anosov diffeomorphism $\varphi:M\to M$ is always partially hyperbolic where the central distribution is the vertical line sub-bundle tangent to the fibers.
\end{fact}
For further discussions and more details see \cite{BuW}. Here is an idea to prove this fact (this is known as the invariant cone method): let $E_u$ be the unstable line sub-bundle on $M$ and define $F_u=\pi^{-1}(E_u)$ (that is, $F_u$ is the inverse image of $E_u$ by the projection). So, $F_u$ is a $\varphi_{\rho}-$invariant plane sub-bundle of $T(M\times \S^1)$ containing the central line sub-bundle $\widetilde{E}_c$ tangent to the fibers. The aim is to construct a $\varphi_{\rho}-$invariant line sub-bundle $\widetilde{E}_u$ such that $\widetilde{E}_c\oplus\widetilde{E}_u=F_u$. For this, we start with any line sub-bundle $E$ in $F_u$ supplementary to $\widetilde{E}_c$ (that is, $\widetilde{E}_c\oplus E=F_u$) and prove that the images of $E$ by action of the iterates of $\varphi_{\rho}$ converge to an invariant line sub-bundle. The way to do this is to show that the action of $\varphi_{\rho}$ on the space of continuous line sub-bundles $E$ in $F_u$ and supplementary to $\widetilde{E}_c$ is in some sense a "contraction" and, therefore, it has a unique fixed point and every orbit converges to it. Constructing $\widetilde{E}_s$ is similar by using $\varphi_{\rho}^{-1}$ instead and we have that $\widetilde{E}_u$ and $\widetilde{E}_s$ verify the desired fact.
\begin{corollary} \label{6.2}
Let $A:\T^2\to \T^2$ be a linear Anosov automorphism and $\rho:\T^2\to \S^1$ be a smooth map. Let $A_{\rho}:\T^2\times\S^1\to \T^2\times\S^1$ be the circle extension of $A$. Then $A_{\rho}$ is partially hyperbolic autonomous diffeomorphism (with respect to a framing).

\end{corollary}
Proof: We have seen (in the previous fact) that $A_{\rho}$ is  partially hyperbolic with central direction tangent to the fibers, and, the stable and the unstable directions define a $\varphi_{\rho}-$invariant horizontal plane field $\mathcal{H}$. Each plane $\mathcal{H}(x)$ projects isomorphically onto $T_{\pi(x)}\T^2$. Let $\mathcal{F}$ be a framing of $\T^2$ for which $A$ is autonomous, then its pullback $\pi_*(\mathcal{F})$ defines two independent vector fields on $\T^2\times\S^1$ contained in $\mathcal{H}$. Let $X$ be a "constant" non-null vector fields on $\T^2\times \S^1$ tangent to the fibers. Then $X$ is preserved by $A_{\rho}$ and $X$ together with $\pi_*(\mathcal{F})$ define a framing of $\T^2\times \S^1$ for which $A_{\rho}$ is autonomous with cocycle matrix equal to $\begin{pmatrix}
1 &0 &0 \\ 
0 &a &b \\
0 &c &d
\end{pmatrix}$ where $\begin{pmatrix}
a &b \\
c &d
\end{pmatrix}$ is the cocycle matrix of $A$ with respect to $\mathcal{F}$.\begin{flushright}
				$\square$
			\end{flushright}
Since the map $\rho:\T^2\to\S^1$ can be anything (for example constant inside some open set and non-constant outside) we deduce that the partially hyperbolic dynamics $A_{\rho}:\T^2\times \S^1\to \T^2\times\S^1$ is not algebraic. Thus, the constructed framing for which $A_{\rho}$ is autonomous cannot be $C^1$ because this would contradict Theorem \ref{3D}. In fact, the lack of regularity occurs in the horizontal plane field generated by the stable and unstable directions.

\section{Appendix: Palais' theorem in the $C^1-$regularity case} \label{appendix}

Here one observes that Palais's Theorem  \ref{palais}   holds in the $C^1$ regularity case. Indeed, its proof is based upon Frobenius' theorem which is true in the $C^1$ case (in fact even in the Lipschitz case (see \cite{sim} and \cite{RAM})). We sketch its proof here and indicate how it actually extends to the $C^1$ setting.\\\\
Let $M$ be a smooth $n-$dimensional manifold and $\Gamma^{\infty}(M)$ be the infinite dimensional algebra of smooth vector fields tangent to $M$. Let $\g$ be a Lie algebra whose associated simply connected Lie group is denoted by $G$. Consider a homomorphism of Lie algebras $X\in\g \to X_M\in\Gamma^{\infty}(M)$ whose image consists of complete vector fields. \\
Denote by $X_G$  the right invariant vector field generated by $X\in \g$ and define the vector field $X_{G\times M}$ on $G\times M$ by $X_{G\times M}(g,x)=(X_G(g), X_M(x))$. 
Let $\phi_{X_{G}}$, $\phi_{X_{M}}$, and $\phi_{X_{G\times M}}$ denote the time-one map flows on $G$, $M$, and $G\times M$ respectively.\\

Define the (local) "pre $G-$action" as follows: for small $h\in G$ we have $\exp^{-1}(h)=X^h\in \g$ and $h$ acts on $G$, $M$, and $G\times M$ as $\phi_{X_{G}^h}$, $\phi_{X_{M}^h}$, and $\phi_{X_{G\times M}^h}$ respectively. Consider the projections on the factors $p_1:G\times M\to G$ and $p_2:G\times M\to M$. 
Then $p_1$ and $p_2$ are equivariant with respect to these pre-actions, that is, $p_i(h(g,x))=hp_i(g,x)$ for $i=1,2$.

 We want to prove that these pre-actions define, indeed,  genuine actions of $G$. The point is that on $G$, we get nothing but the left action of $G$ on itself (Indeed, the flow of a right invariant vector fields acts on $G$ by left multiplications of a one-parameter group). More precisely, $\phi_{X_{G}^h}$ acts as $g\to hg$.\\
Let now $D$ be the distribution on $G\times M$ generated by all vector fields $X_{G\times M}$ for $X\in \g$. Its involutivity follows exactly from the fact that the infinitesimal action $X\to X_M$ is bracket preserving: $[X_{G\times M}, Y_{G\times M}]=([X_G,Y_G], [X_M,Y_M])=([X,Y]_G, [X,Y]_M)$.\\
Let $\mathcal{D}$ be a leaf through some point $(1, x_0)$,  then the projection $p_1: \mathcal D \to G$ is a local diffeomorphism and it conjugates the pre-action on the leaf to the action on $G$. Hence the pre-action on the leaf is in fact a true action. Knowing this,  we use  the projection $p_2: \mathcal D \to M$ to  deduce that the pre-action on $M$ is a true action.

Observe that all the proof can be adapted to the case where the vector fields $X_M$ are not necessarily complete allowing one to get a local action on $M$.\\

Now, an infinitesimal $C^1$ action means that $X_M$ is a $C^1$ vector field and that $[X_M,Y_M]=[X,Y]_M$. All the previous proof extends to this case and yields a (local) action of $G$ on $M$. Indeed, Frobenius' Theorem, that is, involutivity implies integrability i.e. existence of leaves, is valid for $C^1$ distributions (actually  for Lipschitz ones \cite{sim}). Observe that the so-obtained action is $C^1$. To  see this,  remember the classical fact that the local flow of a $C^1$ vector field is $C^1$.

\end{document}